\documentclass[11pt]{amsart}
\usepackage{euscript,amssymb,amsbsy}
\usepackage[all]{xy}
\usepackage{graphicx,bm,a4wide}

\numberwithin{equation}{section}

\newtheorem{stuff}{Stuff}[section]

\newtheorem{theorem}[stuff]{\bf Theorem}
\newtheorem{proposition}[stuff]{\bf Proposition}
\newtheorem{lemma}[stuff]{\bf Lemma}
\newtheorem{corollary}[stuff]{\bf Corollary}

\newenvironment{definition}{%
\vskip1ex\refstepcounter{stuff}\trivlist \itemindent 0pt
\item[\hskip\labelsep\bf Definition \thestuff.]%
\ignorespaces}{\endtrivlist\vskip1ex}%
\newenvironment{example}{%
\vskip1ex\refstepcounter{stuff}\trivlist \itemindent 0pt
\item[\hskip\labelsep\bf Example \thestuff.]%
\ignorespaces}{\endtrivlist\vskip1ex}%
\newenvironment{remark}{%
\vskip1ex\refstepcounter{stuff}\trivlist \itemindent 0pt
\item[\hskip\labelsep\bf Remark \thestuff.]%
\ignorespaces}{\endtrivlist\vskip1ex}%

\let\rar\rightarrow
\let\lar\longrightarrow

\let\mt\mapsto
\let\lmt\longmapsto
\font\tenmsa=msam10 %
\newcommand\hdashpiece{%
{\vrule height2.75pt depth-2.35pt width2.3pt \kern1.7pt}}%
\newcommand\hdashpieces{%
{\hdashpiece\hdashpiece\hdashpiece\hdashpiece}}%
\newcommand\dashto{\mathrel{%
\hdashpiece\hdashpiece\kern-0.4pt\hbox{\tenmsa K}}}%
\newcommand\dashar{\mathrel{%
\hdashpieces\kern-0.4pt\hbox{\tenmsa K}}}%

\let\euf\EuScript %use of package ``euscript'' required
\let\cal\mathcal
\let\mbb\mathbb

\let\mfrak\mathfrak
\DeclareFontFamily{OT1}{rsfs}{}
\DeclareFontShape{OT1}{rsfs}{n}{it}{<->rsfs10}{}
\DeclareMathAlphabet{\crl}{OT1}{rsfs}{n}{it}

\let\unl\underline
\let\unbar\underbar

\let\wtld\widetilde

\let\nit\noindent

\let\disp\displaystyle
\let\srel\stackrel
\let\vphi\varphi

\let\veps\varepsilon

\newcommand\Hom{\mathop{\rm Hom}\nolimits}
\newcommand\End{\mathop{\rm End}\nolimits}
\newcommand\Mor{\mathop{\rm Mor}\nolimits}
\newcommand\Pic{\mathop{\rm Pic}\nolimits}
\newcommand\Spec{\mathop{\rm Spec}\nolimits}

\newcommand\Ker{{\rm Ker}}

\newcommand\invq{{\slash\kern-0.65ex\slash}}
\newcommand\rk{{\rm rank\,}}

\let\l\lambda
\let\L\Lambda

\let\O\Omega
\let\og\omega
\let\si\sigma

\let\th\theta
\let\sm\setminus

\newcommand\s{{\rm s}}
\newcommand\sst{{\rm ss}}
\newcommand\us{{\rm us}}

\newcommand\bbV{{\mbb V}}
\newcommand\bone{{1\kern-0.57ex\rm l}}

\newcommand\cE{{\cal E}}

\newcommand\cO{{\cal O}}

\newcommand\Lie{\mathop{{\cal L}ie}\nolimits}
\newcommand\T{\mathop{\sf T\kern-.2ex}\nolimits}
\newcommand\N{\mathop{\sf N\kern-.1ex}\nolimits}

\newcommand\pr{\mathop{\rm pr}\nolimits}
\newcommand\Sym{\mathop{\rm Sym}\nolimits}
\newcommand\Grs{{\rm Grass}}
\newcommand\codim{{\rm codim}}
\newcommand\Gl{{\rm Gl}}

\newcommand\Ad{{\rm Ad}}

\renewcommand\det{{\rm det}}

\let\ges\geqslant
\let\les\leqslant
\newcommand\ac{{\rm ac}}

\newcommand\btimes{\hbox{\Large$\times$}}
\newcommand\mx{{\rm max}}
\newcommand\mn{{\rm min}}
\newcommand\Flag{{\rm Flag}}
\newcommand\VB{{\cal V\kern-.41ex\euf B}}
\newcommand\bbl{{\mfrak l}}

\def\overset#1#2{\mathop{#2}\limits^{#1}}
\def\underset#1#2{\mathop{#2}\limits_{#1}}
\def\ouset#1#2#3{\underset{#1}{\overset{#2}{#3}}}
\begin{document}

\title{Strong exceptional sequences of vector bundles on certain Fano varieties}

\author{Mihai Halic}


\begin{abstract}
Exceptional sequences of vector bundles over a variety $X$ are special 
generators of the triangulated category $D^b(Coh\,X)$. Kapranov proved 
the existence of tilting bundles over homogeneous varieties for the 
general linear group. King conjectured the existence of tilting sequences 
of vector bundles on projective varieties which are obtained as quotients 
of Zariski open subsets of affine spaces. 

The goal of this paper is to give further examples of strong exceptional 
sequences of vector bundles on certain projective varieties. These are 
obtained as geometric invariant quotients of affine spaces by linear 
actions of reductive groups, as appears in King's conjecture.
\end{abstract}

\maketitle
\markboth{\sc Mihai Halic}
{\sc Strong exceptional sequences of vector bundles}

%%===================================================================

\section*{Introduction}

The concept of derived categories has been introduced by 
Grothen\-dieck and developed further by Verdier. However, their work 
remained within a very general and abstract setting, and people wished 
to have concrete examples which arise from geometry. In algebraic geometry 
one of the essential objects associated to a projective variety is 
the (bounded) derived category of coherent sheaves over it. Its 
knowledge allows to recover all the cohomological data of the variety. 

Beilinson made the first major step by proving that the line bundles 
$
\mathcal O_{\mathbb P^n},\mathcal O_{\mathbb P^n}(1),\!...,
\mathcal O_{\mathbb P^n}(n)
$ 
generate $D^b({\rm Coh}\,\mathbb P^n)$, 
and actually form a tilting sequence. 
Afterwards have appeared several other examples of varieties 
admitting (strong and complete) exceptional sequences of vector 
bundles. One of the most notable results in this direction 
has been obtained by Kapranov \cite{ka}. He explicitly constructed 
tilting sequences of vector bundles over homogeneous varieties 
for ${\rm Gl}(n)$, that is over Grassmannians and flag manifolds.  
Further examples, which are based on Kapranov's result, have been 
obtained in \cite{cm}. 

In the unpublished preprint \cite{ki}, King conjectured that there 
are tilting bundles over projective varieties which are obtained 
as invariant quotients of affine spaces for linear actions of 
reductive groups. Observe that flag varieties for 
$Gl(n,\mathbb C)$, and toric varieties are special cases of such 
quotient varieties. 

The answer to King's conjecture is negative in general. Hille and 
Perling gave in \cite{hp} an example of a toric variety ($\mbb P^2$ 
blown-up successively three times) with the property that it does not 
admit a tilting object formed by line bundles. 
However it is still a very interesting problem to find classes of 
examples for which the conjecture holds. 
In the paper \cite{ah}, Altmann and Hille proved the existence of 
(partial) strong exceptional sequences on toric varieties arising 
from thin representations of quivers, but their construction gives 
sequences of very short length. 

The goal of this paper is to give further examples of strong 
exceptional sequences of vector bundles over certain Fano varieties. 
The varieties considered in this paper are obtained as geometric quotients 
of open subsets of affine spaces by linear actions of a reductive groups. 

For the comfort of the reader, we recall that a sequence of vector 
bundles $(\cal F_1,\ldots,\cal F_z)$ over a variety 
$Y$ is called {\it strongly exceptional} if the following two 
conditions are fulfilled:
\begin{enumerate}
\item $H^0\bigl(Y,\Hom(\cal F_j,\cal F_i)\bigr)=0,
\;\forall\,1\les i< j\les z$; 
\item $H^q\bigl(Y,\Hom(\cal F_j,\cal F_i)\bigr)=0,
\;\forall\,i,j=1,\ldots,z$, and $\forall\,q>0$.
\item A {\em tilting sequence} is a strongly exceptional sequence 
$(\cal F_1,\ldots,\cal F_z)$ with the property that  
$\cal F_1,\ldots,\cal F_z$ generate $D^b({\rm Coh}\,Y)$.
\end{enumerate}
Consider an algebraically closed field $K$ of characteristic zero, 
a connected, reductive group $G$ over $K$, 
and a representation $\rho:G\rar\Gl(V)$. Let 
$\bbV:=\Spec\bigl(\Sym^\bullet V^\vee\bigr)$ be the affine space 
corresponding to $V$. We denote $\chi_\ac=\chi_\ac(G,V)$ the weight 
of the $G$-module $\det V$. We make the following assumptions: 
\begin{enumerate}
\item the ring of invariants $K[\bbV]^T\!=\!K$, 
where $T$ is the maximal torus of $G$;
\item ${\rm codim}_{\bbV}\bbV^\us(G,\chi_\ac)\ges 2$, and 
$G$ acts freely on the semi-stable locus $\bbV^\sst(G,\chi_\ac)$.
\end{enumerate}
\nit We denote $Y:=\bbV\invq_{\chi_\ac}G$ the invariant quotient. 
The main ingredient that we use for constructing exceptional sequences 
over $Y$ is the set ${\cal E_1,\ldots,\cal E_N}$ of `extremal' nef 
vector bundles over $Y$ (see section \ref{sct:nef-vb}). 
They enjoy good cohomology vanishing properties which are required 
by the definition of exceptional sequences. 
The first main result of this paper is the following:\medskip

\nit{\bf Theorem~A}{\kern1ex\it 
Let $d_j:=\rk(\cal E_j)$, and write 
$\chi_\ac=
\mbox{$\overset{N}{\underset{j=1}\sum}\;$} 
m_j\!\cdot\det(\cal E_j)$, with $m_j\ges 1$. 
We denote by $\cal Y_{m,d}$ the set of Young diagrams with at most 
$d$ rows and $m$ columns, and consider  
$$
\begin{array}{ll}
\!\euf{ES}(Y):=
&
\text{the set of vector bundles occurring 
as direct summands in}
\\[.5ex] 
&
\text{the Schur powers } 
\mbb S^{\l^{(1)}}\cal E_1\!
\otimes\ldots\otimes
\mbb S^{\l^{(N)}}\cal E_N, 
\,\l^{(j)}\!\in\cal Y_{m_j-d_j,d_j},
\\[.5ex]
&
\text{corresponding to irreducible $G$-modules}.
\end{array}
$$
Then the vector bundles $\cal E\rar Y\!$, $\,\cal E\in\euf{ES}(Y)$, form 
a strong exceptional sequence over $Y$ with respect to an appropriate 
order on $\euf{ES}(Y)$.

Moreover, if the multiplicities of the isotypical components of $V$ 
are sufficiently high, these vector bundles are slope semi-stable 
with respect to any polarization on $Y$.
}\medskip

The estimates appearing in this theorem are not strong enough to recover 
Kapranov's construction for partial flag varieties. We have to go on, and 
exploit the fibre bundle structure. 
The optimal result would be the following: 

\begin{itemize}
\item[] 
Consider a fibre bundle $Y\srel{\phi}{\rar}X$. 
Suppose that $(\cal F_i)_{i\in I}$ is a strong 
exceptional sequence of vector bundles on $X$, and that 
$(\cal E_j)_{j\in J}$ is a sequence of vector bundles on $Y$ 
whose restriction to the fibres of $\phi$ give rise to strong 
exceptional sequences relative to $\phi$. 
Then $(\phi^*\cal F_i\otimes\cal E_j)_{(i,j)\in I\times J}$ 
is a strong exceptional sequence on $Y$. 
\end{itemize}

\nit Unfortunately such a statement is overoptimistic in general. 
The content of our second main result is that the statement 
above becomes true under suitable restrictive hypotheses on 
the fibration $\phi$. 
More precisely, we place ourselves in the following framework:
\begin{enumerate}
\item[(i$'$)] There is a quotient group $H$ of $G$ with kernel $G_0$, 
and a quotient $H$-module $W$ of $V$ 
with kernel $V_0$, such that the natural projection 
$\pr^{\bbV}_{\mbb W}:\bbV\rar\mbb W$ has the following property: 
 
\centerline{
$\pr^{\bbV}_{\mbb W}\bigl(\;\bbV^\sst\bigl(G,\chi_\ac(G,V)\bigr)\;\bigr)
\subseteq\mbb W^\sst\bigl(G,{\chi_\ac}(H,W)\bigr).$
} 
\nit We denote by $Y\srel{\phi}{\rar}X$ the induced morphism 
at the quotient level. 

\item[(ii$'$)] The unstable loci have codimension at least two, and both 
quotients\smallskip 

\centerline{
$\bbV^\sst(G,\chi_\ac(G,V))\rar Y\;$ and 
$\;\mbb W^\sst(H,\chi_\ac(H,W))\rar X$} 
\nit are principal bundles.

\item[(iii$'$)] The nef cone of the total space $Y$ is the sum of the nef cones 
of the base $X$, and that of the fibre: 
$\crl N(G,V)=\crl N(H,W)+\crl N(G_0,V_0)$. Denote $\VB^+(X)$ and 
$\VB^+_0$ the corresponding sets of extremal nef vector bundles. 

\item[(iv$'$)] The maximal torus $T_0\subset G_0$ has exactly $\dim T_0$ 
weights on $V_0$.
\end{enumerate}

\nit Our main result in the relative case is the following:\medskip 

\nit{\bf Theorem~B}
{\kern1ex\it Let us denote $d_F\!:=\!\dim F,$ for $F\in\VB^+(X)$, 
and $d_E\!:=\!\dim E,$ for $E\in\VB^+_0$. We write 
$
\chi_\ac(H,W)=\!
\mbox{$\underset{F\in\VB^+(X)}\sum$}\kern-1.5ex
m_F\cdot\det F\;(m_F\ges 0),
$ 
and  
$\chi_\ac(G_0,V_0)=\!
\mbox{$\underset{E\in\VB^+_0}\sum$} 
m_E\cdot\det E\;\;(m_E\ges 0).
$ 
Suppose $(b_F)_{F\in\VB^+(X)}$ are integers such that for all $q>0$, 
and for all Young diagrams $\beta^{F}$ of length $d_F$, with 
$\beta^{F}_\mn\ges -b_F$,  holds: 
$
H^q\biggl(
X,\mbox{$\underset{F\in\VB^+(X)}\bigotimes$}\kern-1ex 
\mbb S^{\beta^{F}}\cal F
\biggr)=0.
$ 

Then the elements of the set ${\euf{ES}}(Y)$ defined below 
form a strong exceptional sequence of vector bundles over $Y$: 
$$
\begin{array}{ll}
{\euf{ES}}(Y):=
&
\text{all the direct summands, corresponding to irreducible}
\\[0ex] 
&
\text{$G$-modules contained in }
\\[1ex] 
&
\phi^*\bigl(\mbb S^{\l^{\bullet}}\cal F_\bullet\bigr)
\otimes
\mbb S^{\nu^{\bullet}}\cal E_\bullet
:=
\phi^*\bigl(
\mbox{\kern-2ex$\underset{\tiny F\in\VB^+(X)}\bigotimes$}
\kern-1.7ex\mbb S^{\l^{F}}\cal F\,
\bigr)
\otimes
\mbox{$\underset{E\in\VB^+_0}\bigotimes$\,}
\kern-.5ex\mbb S^{\nu^{E}}\cal E, 
\end{array}
$$
with $\l^F\in\cal Y_{b_F,\,d_F}$ and $\nu^{E}\in\cal Y_{m_E-d_E,\,d_E}.$ 

\nit Moreover, it holds: 
$H^q\biggl(
Y,\underset{F\in\VB^+(X)}\bigotimes\kern-1.7ex
\phi^*\,\mbb S^{\beta^{F}}\cal F
\;\otimes
\underset{E\in\VB^+_0}\bigotimes\kern-.5ex
\mbb S^{\alpha^{E}}\cal E
\biggr)=0$ 
for all $q>0$, and all Young diagrams $\beta^{F}$ and $\alpha^{E}$ 
of length $d_F$ and $d_E$, 
with $\beta^{F}_\mn\ges -b_F$ and $\alpha^{E}_\mn\ges -(m_E-d_E)$ 
respectively.
}\medskip

We point out that in both cases it {\em remains open} the question 
under which hypothesis these sequences are/extend to {\em tilting} 
objects. However, we remark that, taking into account the 
example constructed in \cite{hp}, a {\em general answer} concerning 
the (non-)existence of tilting vector bundles over quotients of 
affine spaces must be involved.

The definition of an exceptional set involves two conditions. 
Accordingly, the paper is divided in two main parts, each focusing 
on one of the two conditions: 

-- The sections \ref{sct:stability} and \ref{sct:conseq-stab} form the 
first part: we prove a stability result for associated vector bundles, 
and define an order on the set of irreducible $G$-modules 
for which there are no homomorphisms from a `larger' vector bundle 
into a `smaller' one (see theorem \ref{thm:h000}). 

-- The sections \ref{numer-crit}, \ref{sct:nef-vb} and \ref{cohom-nef}, 
have a preparatory character: we introduce the `extremal' nef vector 
bundles, and study their cohomological properties. 

-- The second part of the article consists of the sections 
\ref{sct:main} and \ref{sct:main2}: they contain the proofs 
of the main results. 
The main tool used for proving the vanishing of the higher cohomology groups 
is a result due to Manivel (see \cite{ma}), and Arapura (see \cite{ar}). 
However, this general result is not sufficient to address the relative case, 
and we have to dwell on our particular context. 
In theorem \ref{thm:direct-image} we prove the following nefness property, 
which is an essential ingredient in the proof of Theorem B.\medskip 

\nit{\bf Theorem} 
{\kern1ex\it 
Suppose that $\mbb V\rar\mbb W$ satisfies the properties 
\mbox{\rm(i$'$)} and \mbox{\rm(ii$'$)} above, and denote 
$Y\srel{\phi}{\rar}X$ the morphism induced at the quotient level. 

Let $\cal E\rar Y$ be a nef vector bundle, associated to a $G$-module $E$. 
Then $R^q\phi_*\cal E=0$ for all $q>0$, 
and $\phi_*\cal E\rar X$ is still a nef vector bundle. 
}\medskip

-- Finally, in section \ref{sct:expl}, we illustrate the general theory. 
On one hand, we recover Kapranov's construction for the Grassmannian and 
for flag varieties, by using our results. 
On the other hand, we give further examples of strong exceptional 
sequences over quiver varieties. The very pleasant feature is that 
we obtain these example by an almost algorithmic procedure, which 
applies to any quiver variety. 

Some of the results have been presented at the HOCAT 2008 Conference, 
held at Centre de Recerca Matem\`atica, Bellaterra, Spain. 
\medskip

%\nit{\sl Acknowledgment}\quad The author wishes to acknowledge KFUPM for its support ++++++ 

%%%%%%%%%%%%%%%%%%%%%%%%%%%%%%%%%%%%%%%%%%%%%%%%%%%%%%%%%%%%%%%%%%%%%

\section{A stability property}
{\label{sct:stability}}\setcounter{equation}{0}

The symbol $\mbb Q$ will always denote the field of rational numbers, 
and $K$ will be an algebraically closed field $K$ of characteristic zero. 
Throughout the paper, $G$ will always denote a connected, reductive 
group over $K$, and $T$ will be the maximal torus of $G$. 
We consider a faithful representation $\rho:G\rar\Gl(V)$, and denote 
by $\bbV:=\Spec(\Sym^\bullet V^\vee)$ the corresponding affine space. 
We shall assume that the ring of invariants $K[\bbV]^T=K$; 
it follows automatically that $K[\bbV]^G=K$. 

\begin{lemma}{\label{lm:Zneq0}}
Let $V$ be a non-zero $G$-module such that $K[\bbV]^T=K$. Then: 
\begin{enumerate} 
\item There is a 1-PS $\l\in\cal X_*(T)$ such that all 
its weights on $V$ are strictly positive.  
\item $G$ is not semi-simple. 
\end{enumerate}
\nit We fix once for all 
$\bbl\in\cal X_*(T)\otimes_{\mbb Z}\mbb R$ 
such that its weights on $V$ are all positive, and moreover it has 
`irrational slope', that is $\Ker(\bbl:\cal X^*(T)\rar\mbb R)=\{0\}.$ 
\end{lemma}

\begin{proof}
(i) Let $\Phi$ denote the set of weights of the $T$-module $V$. 
Then the set of weights of the $T$ on $K[\bbV]$ is the `cone' 
$\underset{\eta\in\Phi}\sum\mbb N\eta$. 
Since $K[\bbV]^T=K$, this cone is strictly convex. 
Otherwise we can construct a non-trivial $T$-invariant monomial. 
It follows that there is $\l\in\cal X_*(T)$ with 
$\langle\eta,\l\rangle>0$ for all $\eta\in\Phi$.\smallskip 

\nit (ii) Assume that $G$ is a semi-simple group. 
The previous step implies that $K[\bbV^m]^T=K$, 
hence $K[\bbV^m]^G=K$ for all $m\ges 1$. 
Since $G$ is semi-simple, it has an 
open orbit in $\bbV^m$. For large $m$ we get a contradiction. 
\end{proof}

Let $\th\in\cal X^*(G)$ be a character. We denote: 
\begin{equation}{\label{G-sst}}
\begin{array}{l}
K[\bbV]^G_\th:=\{f\in K[\bbV]\mid f(g\times y)=\th(g)\cdot f(y),
\,\forall y\in\bbV\}
\\[2ex]
K[\bbV]^{G,\th}:=K\oplus\underset{n\ges 1}\bigoplus K[\bbV]^G_{\th^n},
\\[2ex] 
\bbV^\sst(G,\th):=\{y\in\bbV\mid\exists n\ges 1\text{ and }
f\in K[\bbV]^G_{\th^n}\text{ s.t. }f(y)\neq 0\}.
\end{array}
\end{equation}
We say that $\th$ is {\em effective} if there is $n\ges 1$ such that 
$K[\bbV]^G_{\th^n}\neq 0$, that is $\bbV^\sst(G,\th)\neq\emptyset$. 

\begin{definition}
We define the {\it anti-canonical character} of the $G$-module $V$ to 
be the character of the $G$-module $\det V$. 
\end{definition}
Explicitly: decompose 
$V=\underset{\og\in\cal X}\bigoplus M_\omega^{\oplus m_\og}$ 
into its $G$-isotypical components. Let $\chi_\og$ be the character 
by which $Z(G)^\circ$ acts on $M_\og$, and denote $d_\og:=\dim M_\og$. 
Then 
$\;
\chi_\ac(G,V):=
\mbox{$\underset{\og\in\cal X}\sum$}
m_\og  d_\og\chi_\og\in\cal X^*(G). 
$ 
For shorthand, we will write $\chi_\ac=\chi_\ac(G,V)$. 

\begin{lemma}{\label{lm:effective}}
Assume that $m_\og \ges d_\og$. Then the character $\chi_\og$ is effective. 
Moreover, if $m_\og > d_\og$ for all $\og$, then $\chi_\ac$ is 
effective, and the $\chi_\ac$-unstable locus has codimension at least two.
\end{lemma}

\begin{proof}
We view $V$ as $\bigoplus_{\og\in\cal X} \Hom(K^{m_\og }, M_\og)$. 
Since $m_\og \ges d_\og$, we can associate to an element 
$\Hom(K^{m_\og }, M_\og)$ the $d_\og\times d_\og$-minor 
corresponding to the first $d_\og$ columns. This defines a regular 
function $f_\og$ which is $d_\og\chi_\og$-equivariant; moreover, $f_\og$ 
does not vanish on surjective homomorphisms. 
It follows that $d_\og\chi_\og$, and therefore $\chi_\og$, is effective 
for all $\og$. 

If a point belongs to the unstable locus, then all the minors $f_\og$ 
have to vanish. Since $m_\og \ges d_\og+1$, this implies the vanishing 
of at least two independent minors.
\end{proof}

Now we prove a general stability result of independent interest. 
It is well known that the tangent bundle of the projective space 
is stable, and more generally the tautological bundles over 
Grassmannians are stable. Our goal is to generalize these facts. 

We denote ${\{G_j\}}_{j\in J}$ the simple factors of $G$, and let 
$\gamma_j:G\rar G_j$, %for $j\in J$, 
be the corresponding quotient morphisms. 
Using the $\gamma_j$'s we extend the structural group of $\O\rar Y$, and 
obtain the principal $G_j$-bundles $\O(G_j)\rar Y$. 
The main result of this section is:

\begin{theorem}{\label{thm:stab-bdl}}
Assume that $G$ acts freely on $\O\!:=\!\bbV^\sst(G,\th)$, for some 
$\th\in\cal X^*(G)$, and let $Y$ be the quotient. 
Assume that $m_\omega >\dim M_\og$ holds for all $\omega\in\cal X$. 
Then the principal $G_j$-bundles $\O(G_j)\rar Y\!$, $j\in J$, obtained 
by extending the structural group are semi-stable. 
\end{theorem}

\begin{proof}
We fix $j\in J$, and a maximal parabolic subgroup $P_j\subset G_j$; denote 
$P:=\gamma_j^{-1}P_j$: it is a maximal parabolic subgroup of $G$. 
We observe 
that the associated homogeneous bundles $\bigl(\O(G_j)\bigr)(G_j/P_j)$ 
and $\O\bigl(G/P\bigr)$ are isomorphic. 

We denote $H=\underset{\omega }{\prod}\,H_\omega :=
\underset{\og }{\prod}\,\Gl_K(m_\omega )$: it acts naturally on 
$\bbV$; the $G$- and $H$-actions on $\bbV$ commute. It follows 
that $H$ still acts on $\O(G/P)$ by  

\centerline{$
H\times\O(G/P)\rar\O(G/P),\; h\times[y,gP]:=[hy,gP].
$}

\nit We will prove that whenever there is a reduction of the 
structural group 

\centerline{$
s:Y^o\rar\bigl(\O(G_j)\bigr)(G_j/P_j)=\O\bigl(G/P\bigr),
$ 
with $Y^o\subset Y$ open and $\codim_Y(Y\sm Y^o)\ges 2,$}

\nit holds $\deg_{Y}\big(s^*{\sf T}_{\O(G/P)/Y}\big)\ges 0.$ 
Equivalently, the reduction $s$ can be viewed as a 
$G$-equivariant morphism $S:\O^o=q^{-1}(Y^o)\rar G/P$. 

The idea is to move $s$ using the $H$-action on $\O(G/P)$. 
Let $\hat y\in Y$ be a generic point, and consider $y\in\O$ 
over $\hat y$. We define the following subgroups of $H$: 
$\quad K_{\hat y}:={\rm Stab}_H(y)$, and 
$$ 
H_{\hat y}:=\{h\in H\mid \exists\,g_h\in G\text{ s.t. }
h\times y=\rho(g_h^{-1})y\}
=\mbox{$\underset{\og}\prod$}  H_{\og ,\hat y}.
$$
We observe that $K_{\hat y}$ does not depend on the choice of 
$y\in q^{-1}(\hat y)$. 
Since $G$ acts freely on $\O$, the assignment $h\mt g_h$   
defines a group homomorphism $\rho_{\hat y}:H_{\hat y}\rar G$ whose 
kernel is $K_{\hat y}$. We move the section $s$ using the action 
of $H_{\hat y}$. For $h\in H_{\hat y}$ define a new section 
$s_h$ as follows: 
$$
s_h(\hat x):=[x,S(h^{-1}\times x)]
\quad
\text{(equivalently, $S_h(x):=S(h^{-1}\times x)$).}
$$
Observe that as $h\in H_{\hat y}$ varies, 
$s_h(\hat y)=h\times s(\hat y)$ moves in the vertical direction.
\medskip 

\nit\unl{Claim}\quad $H_{\hat y}/K_{\hat y}\rar G/Z(G)^\circ$ is 
surjective. Write $y={(y_\omega )}_\omega $ w.r.t. the 
direct sum decomposition of $V$; for each $\omega\in\cal X$, 
$y_\omega =(y_{\omega 1},\ldots,y_{\omega m_\omega})$. Since $y\in\O$ is 
chosen generically, and $m_\omega>\dim M_\og=:d_\omega$, we may 
assume that for each $\omega\in\cal X$ the vectors 
$y_{\omega 1},\ldots,y_{\omega m_\omega}$ span $M_\og$. Equivalently, 
we may view $y_\omega$ as a surjective homomorphism 
$K^{m_\omega}\rar M_\og$. 

For $g\in G$ holds $\rho(g)y={(\rho_\omega (g)y_\omega)}_\omega$. Using 
that $m_\omega> d_\omega$, we deduce that for each $\omega\in\cal X$ 
there is $h_\omega\in\Gl_K(m_\omega)$ such that 
$h_\omega y_\omega=\rho_\omega(g^{-1})y_\omega$. 
For $h:={(h_\omega)}_\omega$ we have $hy=\rho(g^{-1})y$, that is 
$g\in{\rm Image}\bigl(H_{\hat y}/K_{\hat y}\rar G\bigr)$. 
\medskip 

Back to our proof: the infinitesimal action of 
$H_{\hat y}$ preserves the restriction to the fibre 
$q^{-1}(\hat y)=\{[y,gP]\mid g\in G\}\cong G/P$ of the relative 
tangent bundle ${\sf T}_{\O(G/P)/Y}$. 
By this isomorphism the relative tangent bundle corresponds to 
${\sf T}_{G/P}\rar G/P$. 
The claim implies that the infinitesimal action 
$\Lie(H_{\hat y})\rar{\sf T}_{\O(G/P)/Y, s(\hat y)}$ is surjective. 
Hence there is a section 
$\si\in H^0(Y^o,s^*\det{\sf T}_{\O(G/P)/Y})$ which does not vanish 
at $\hat y$. It follows 
$\deg_Y\big(s^*{\sf T}_{\O(G/P)/Y}\big)\ges 0$. 
\end{proof}

\begin{corollary}{\label{cor:stab-bdl}}
Assume $\th\in\cal X^*(G)$ has the property that $G$ acts freely 
on $\O:=\bbV^\sst(G,\chi)$, and let $Y$ be the quotient. 
Let $E$ be an irreducible $G$-module, and denote by $\cal E:=\O(E)$ 
the associated vector bundle over $Y$. 
Assume that $m_\omega>\dim M_\og$ holds for all $\omega\in\cal X$. 
Then $\cal E\rar Y$ is slope semi-stable with respect to the 
polarization induced by the character $\th$. 
\end{corollary}

\begin{proof} We may assume that 
$G=Z(G)^\circ\times\bigl(\underset{j\in J}{\times}G_j\bigr)$. Since 
each $\O(G_j)$ is semi-stable, $\O\rar Y$ itself is semi-stable. 
The homomorphism $\rho_\omega:G\rar\Gl(E)$ maps $Z(G)^\circ$ into 
the centre of $\Gl(E)$. By using \cite[theorem 3.18]{rr}, we deduce 
that $\cal E=\O(E)\rar Y$ is semi-stable. 
\end{proof}

%%%%%%%%%%%%%%%%%%%%%%%%%%%%%%%%%%%%%%%%%%%%%%%%%%%%%%%%%%%%%%%%%%%%%

\section{The $H^0$ spaces}
{\label{sct:conseq-stab}}\setcounter{equation}{0}

Assume that $E$ is a $G$-module. We will denote by $\cal E$ the vector 
bundle over $Y$ associated to it. More precisely, $\cal E$ corresponds 
to the module of covariants $\bigl(K[\bbV]\otimes_K E^\vee\bigr)^G$. 

The classical Schur lemma says that for two irreducible $G$-modules 
$E$ and $F$, the space $\Hom(E,F)$ consists either of scalars (if $E=F$), 
or vanishes (if $E\neq F$). 
In this section we will prove that a similar result holds 
for the associated vector bundles $\cal E$ and $\cal F$. 

For warming-up, we start with a special case. 
We have proved in corollary \ref{cor:stab-bdl} that $\cal E\rar Y$ is a 
semi-stable vector bundle w.r.t. any polarization on $Y$, as soon 
as the multiplicities $m_\og> d_\og$ for all $\og$. 
Its first Chern class equals $\dim(E)\cdot\chi_E$, where $\chi_E$ denotes 
the character of $Z(G)^\circ$ on $E$. 
Let $\th\in\cal X^*(G)$ be an ample class on $Y$; the slope of $\cal E$ 
w.r.t. $\th$ equals 
$$ 
\mu_\th(\cal E)=
\frac{\deg_\th\cal E}{\dim E}=
\langle \chi_\og\cdot\th^{\dim Y-1},[Y] \rangle. 
$$

\begin{definition}{\label{defn:order1}}
Let $\th$ be a polarization of $Y$. We define the order $<_\th$ on 
$\cal X^*\bigl(Z(G)^\circ\bigr)$ as follows: 
we declare that $\chi<_\th\eta$ if holds: 
$$
\mu_\th(\chi):=\langle\chi\cdot\th^{\dim Y-1},[Y]\rangle
< 
\mu_\th(\eta):=\langle\eta\cdot\th^{\dim Y-1},[Y]\rangle.
$$
Observe that, by the hard Lefschetz property, we can choose 
$\th$ in such a way that 
$\chi=\eta\Leftrightarrow\mu_\th(\chi)=\mu_\th(\eta)$. 
\end{definition}

\begin{proposition}{\label{prop:h00}}
We assume that $m_\og> d_\og$ holds for all $\og$. 
Let $E$ and $F$ be two distinct irreducible $G$-modules, 
such that $Z(G)^\circ$ acts on them by two different characters 
$\chi_E$ and $\chi_F$ respectively, such that 
$\mu_\th(\cal E)<\mu_\th(\cal F)$. 
Then $H^0\bigl(Y,\Hom(\cal F,\cal E)\bigr)=0$. 
\end{proposition}

\begin{proof}
This is an immediate consequence of the semi-stability property 
of $\cal E$ and $\cal F$. 
\end{proof}

The proposition has two shortcomings: first, we have imposed the 
condition on the multiplicities; second, there are distinct representations 
$E$ and $F$ such that the characters $\chi_E$ and $\chi_F$ coincide. 
So we need to sharpen our result.

\begin{theorem}{\label{not-eff}}
Assume that ${\rm codim}_\bbV\bbV^\us(G,\chi_\ac)\ges 2$. 
Let $E$ be an irreducible $G$-module, and let $\cal E\rar Y$ be the 
associated vector bundle. Suppose that there is a weight $\veps$ of 
$T$ on $E$ which is not $T$-effective 
(that is $\bbV^\sst(T,\veps)=\emptyset$). Then $H^0(Y,\cal E)=0$.
\end{theorem}

\begin{proof}
Recall that $H^0(Y,\cal E)=\Mor(\mbb V\rar E)^G$, where 
$$
(g\times S)(y)=g\times S(g^{-1}\times y),\quad\forall\,g\in G
\text{ and }\bbV\srel{S}{\rar}E.
$$ 
Assume that there is 
a non-zero $G$-equivariant morphism $S:\mbb V\rar E$. Then the linear 
span $\langle S\rangle:=\langle S(y), y\in\mbb V\rangle$ is actually a 
$G$-submodule of $E$. Since $E$ is irreducible and $S\neq 0$, 
we deduce $\langle S\rangle=E$. 

On the other hand, $\veps$ is a weight of $T$ on $E$ which is not 
effective. We choose a one dimensional $T$-submodule $E_\veps\subset E$, 
and consider {\em the function} $S_\veps:=\pr^E_{E_\veps}\circ\, S$. 
Then  
$
S_\veps(t\times y)=\veps(t)\cdot S_\veps(y),\;
\forall t\in T,\,y\in\mbb V.
$

Since $\veps$ is not effective, the function $S_\veps$ must vanish. 
This implies that the image of the morphism $S$, and consequently its 
linear span $\langle S\rangle$, is contained in the complement $E'$ of 
$E_\veps$. The contradiction shows that $\langle S\rangle=E$. 
\end{proof}

In order to check that a sequence of vector bundles forms an exceptional 
sequence, one has to prove that there are no non-trivial homomorphisms 
from `larger' bundles into `smaller' ones. Now we define the total order 
required for this property.

\begin{definition}{\label{defn:order2}}
Consider $\bbl\in\cal X_*(T)$ as in lemma \ref{lm:Zneq0}. 
\begin{enumerate} 
\item For any irreducible $G$-module, we define 
$$\mfrak l(E):=
\max\{\langle\eta,\mfrak l\rangle\mid\eta\text{ is a weight of $T$ on $E$}\}.
$$ 
Equivalently:\\ 
$\mfrak l(E)=\langle\eta_E,\mfrak l\rangle$, 
where $\eta_E$ is the  dominant weight of $E$ 
(with respect to $\mfrak l$).  
\item Let $E$ and $F$ be two irreducible $G$-modules. We say that 
$E<_{\bbl}\,F$ if $\mfrak l(E)\,<\,\mfrak l(F)$. 
\end{enumerate}
\end{definition}
Since $\mfrak l$ has irrational slope, for any two irreducible 
$G$-modules $E$ and $F$ holds: 

\centerline{$\mfrak l(E)=\mfrak l(F)\Rightarrow E=F$.}
\nit Hence $<_{\bbl}$ is a total order relation.

The following result can be viewed as a generalization of Schur's lemma.

\begin{theorem}{\label{thm:h000}}
Assume that ${\rm codim}_\bbV\bbV^\us(G,\chi_\ac)\ges 2$. 
\begin{enumerate} 
\item Let $E$ be an irreducible $G$-module. 
Then $H^0\bigl(Y,\End(\cal E)\bigr)=K.$
\item Let $E$ and $F$ be two irreducible $G$-modules such that 
$E<_{\bbl}\,F$. Then\smallskip 

\centerline{$H^0\bigl(Y,\Hom(\cal F,\cal E)\bigr)=0.$}
\end{enumerate}
\end{theorem}

\begin{proof}
(i) A section $s\!\in\! H^0\bigl(Y,\End(\cal E)\bigr)$ corresponds to a 
$G$-equivariant morphism $S:\bbV\rar\End(E)$, where the action on 
$\End(E)$ is by conjugation. 
(Here we use the hypothesis on the codimension of the 
unstable set: regular maps defined on the semi-stable locus 
extend to the whole affine space.)
We will prove that the morphism $S$ is a 
scalar multiple of the identity. 

The origin $0\in\bbV$ is fixed under $G$. Since $S$ is $G$-equivariant, 
the homomorphism $S_0\in\End(E)$ is $\Ad_G$-invariant. 
Schur's lemma implies that $S_0=c\cdot\bone_{E}$, with $c\in K$. 
By lemma \ref{lm:Zneq0}, there is a 1-PS 
$\l\in\cal X_*(T)$ such that all its weights on $V$ are strictly 
positive. In particular $\disp\lim_{t\rar 0}\l(t)y=0$ for all $y\in\bbV$. 
The $G$-equivariance implies $S_{\l(t)y}=\Ad_{\l(t)}\circ S_y$, hence 
$\,\disp\lim_{t\rar 0}\Ad_{\l(t)}\circ S_y =S_0=c\bone_{E}.$

The $\l(t)$-action on $E$ can be diagonalized in an appropriate 
basis formed by weight vectors. We denote ${\{E_i\}}_{i\in I}$ 
the weight spaces of $E$.  
We order the elements of $I$ in decreasing order, and consider the 
corresponding basis in $E$. Then w.r.t. this basis, $S_y$ 
has the following block-matrix shape:
$$
S_y
=
\left(
\begin{array}{c|c|c}
c\bone & *&*\\ \hline 
0&c\bone&*\\ \hline
0&0&c\bone
\end{array}
\right)
\;\text{or equivalently}\;
S_y-c\bone
=
\left(
\begin{array}{c|c|c}
0 & *&*\\ \hline 
0&0&*\\ \hline 
0&0&0
\end{array}
\right),
\;\forall\,y\in\bbV
$$
Let $\mfrak N_\l$ be the vector space which is formed by matrices 
having this shape ($\mfrak N_\l$ is actually a nilpotent Lie algebra). 
Intrinsically, 
$$
\disp\mfrak N_\l=\{A\in\End(E)\mid\lim_{t\rar 0}\Ad_{\l(t)}\circ A=0\}.
$$ 
We denote 
$
\Ker(\mfrak N_\l):=\kern-1ex
\mbox{$\underset{N\in\mfrak N_\l}\bigcap$}\kern-1ex\Ker(N).
$ 
By Engel's theorem, $\Ker(\mfrak N_\l)$ is a non-zero 
vector subspace of $E$. 
Applying the $G$-equivariance once more, we deduce that for any $g\in G$ 
holds: 
$$
Ad_{g^{-1}}\circ
\bigl(
S_y-c\bone
\bigr)
=S_{g^{-1}y}-c\bone\in\mfrak N_\l.
$$
It follows that for all $g\in G$,
$$
\Ker
\bigl(
S_y-c\bone
\bigr)
\supset
g\cdot\Ker(\mfrak N_\l)
\;\Longrightarrow\;
\Ker
\bigl(
S_y-c\bone
\bigr)
\!\supset\kern-.7ex
\mbox{$\underset{g\in G}\sum$}
g\cdot\Ker(\mfrak N_\l).
$$
Note that the right-hand-side is a non-zero $G$-submodule of $E$. 
Since $E$ is irreducible, it follows taht $\Ker\bigl(S_y-c\bone\bigr)=E$, 
that is $S_y=c\bone$ for all $y\in\bbV$.

\nit (ii) The $G$-module $\Hom(F,E)=F^\vee\otimes E$ contains  
the difference $\veps:=\eta_E-\eta_F$ of the corresponding 
dominant characters. Since $E<_\bbl\,F$, $\bbl(E)-\bbl(F)<0$, 
the weight $\veps$ is not $T$-effective. The conclusion follows 
from theorem \ref{not-eff}.
\end{proof}

%%%%%%%%%%%%%%%%%%%%%%%%%%%%%%%%%%%%%%%%%%%%%%%%%%%%%%%%%%%%%%%%%%%%%

\section{Numerical criteria for semi-stability}
{\label{numer-crit}}\setcounter{equation}{0}

In this section we are reviewing some numerical criteria for 
semi-stability, needed later on. 
The following convention is used throughout this section: 
the letters $E, V, W$ denote $G$-modules, 
while the symbols $\mbb{E, V, W}$ will denote 
the corresponding affine spaces: {\it e.g.} 
$\mbb E:=\Spec\bigl(\Sym^\bullet E^\vee\bigr)$. 

For a $G$-module $W$, let $\eta_1,\ldots,\eta_R$ be the weights of the 
maximal torus $T\subset G$. We define:  
$\;m:\mbb W\times \cal X_*(G)_{\mbb R}\rar \mbb R,$ 
$$
\begin{array}{l}
m(w,\l):=
\min
\biggl\{\;
j\;\biggl|\;
\begin{array}{l}
\text{the $t^j$-isotypical component of 
$w$ w.r.t. $\,\l$}
\\ 
\text{does not vanish}
\end{array}
\biggr.\biggr\}.
\end{array}
$$
Observe that for $\l\in\cal X_*(T)$ holds: 
$$
m(w,\l):=
\min
\biggl\{
\langle\eta_j,\l\rangle
\;\biggl|\;
\begin{array}{l}
\text{the $\eta_j$-isotypical component of $w$}
\\ 
\text{does not vanish}
\end{array}
\biggr\}.
$$
We fix a norm $|\cdot|$ on $\cal X_*(T)$, invariant under the Weyl group 
of $G$. For a character $\th\in\cal X^*(G)$, the Hilbert-Mumford criterion 
for $(G,\th)$-(semi-)stability reads: 
\begin{eqnarray}{\label{fct-m}}
\begin{array}{rl}
w\!\in\!\mbb W^{\s\,{\rm (resp. }\,\sst\rm)}(G,\th)
&\Leftrightarrow
m(w)\!:=\!
\inf
\left\{
\frac{\langle \th,\l\rangle}{|\l|}
\,\biggl|\, m(w,\l)\ges 0\biggr.
\right\}
\underset{(\ges)}{>} 0
\\[2ex]
&\Leftrightarrow 
\biggl[\,
m(w,\l)\ges 0 
\,\Rightarrow\, 
\langle \th,\l\rangle\underset{(\ges)}{>} 0
\,\biggr].
\end{array}
\end{eqnarray} 
For $w\in\mbb W$ we define:  
$$
\begin{array}{rl}
S(w):=&
\{\eta_j\mid
\text{the $\eta_j$-isotypical component of $w$ does not vanish}
\}
\\[1.5ex] 
\crl C_w=&
\underset{\eta\in S(w)}\sum\mbb R_{\ges 0}\eta
\\[1.5ex]
\L^G_w:=&
\{\l\in\cal X_*(G)\mid m(w,\l)\ges 0\} 
\\[1.5ex]
\L^T_w:=&
\{\l\in\cal X_*(T)\mid m(w,\l)\ges 0\}
\\[1.5ex]
=&
\{\l\in\cal X_*(T)\mid 
\langle\eta,\l\rangle\ges 0,\;\forall\eta\in\crl C_w
\}
=\crl C_w^\vee.
\end{array}
$$ 
Note that $\crl C_w$ and $\L^T_w$ are convex, polyhedral cones. 
Since there are finitely many $\eta$'s, only finitely many 
cones $\crl C_w$ and $\L_w^T$ occur as $w$ varies in $\mbb W^s(G,\th)$. 
We are interested in the {\em minimal} cones $\crl C_w$. 

\begin{definition}
Let $\th$ be a character of $G$. 
A subset $S\subset\{\eta_1,\ldots,\eta_R\}$ is 
{\em minimal for $\th$} if 
$
\th\in
\mbox{$\underset{\eta\in S}\sum$}
\mbb R_{\ges 0}\eta
$ 
and 
$
\th\not\in
\mbox{$\underset{\eta\in S\sm\{\eta_0\}}\sum$}
\mbb R_{\ges 0}\eta
$
for all $\eta_0\in S$.
\end{definition}
We denote $S_1,\ldots,S_z$ the (finitely many) minimal sets for $\th$, 
and the corresponding cones by 
$\crl C_j$ and $\L_j:=\crl C_j^\vee$, $j=1,\ldots,z$, respectively. 
The Weyl group of $G$ operates by permutations on them. 

Observe that 
$\L^G_w=
\underset{g\in G}\bigcup\Ad_{g^{-1}}\bigl(\L^T_{gw}\bigr)$. 
As $\th$ is $\Ad_G$-invariant, the numerical criterion can 
be reformulated as follows: 
\begin{equation}{\label{intersect}}
\th\in\cal X^*(G)\cap{\rm int.}
\biggl(\,
\bigcap_{w\in\mbb W^\s(G,\th)}
\kern-2ex\crl C_w
\biggl)
=
{\rm int.}\biggl(
\cal X^*(G)\cap\crl C_1\cap\ldots\cap\crl C_z
\biggr). 
\end{equation} 

For two $G$-modules $V, E$, we define the 
$K^\times\times G$-module $W_m:=E\times V^{m}$, $m\ges 1$, 
with the module structure given by

\centerline{$
(t,g)\times\bigl(\vphi,(v_j)_j\bigr)
:=
\bigl(t\cdot(g\times\vphi),(g\times v_j)_j\bigr). 
$}

\nit Consider $l>0$, and define 
$\th_m:=l\chi_t+m\chi_\ac\in\cal X^*(K^\times\times G)$. 
The numerical functions on $\mbb V^m$ and $\mbb W_m$ are the following:
$$
\begin{array}{l}
\disp 
m(\unl v,\l)
\,=\,
\min_j m(v_j,\l),\quad \forall\,\unl v=(v_j)_j\in\mbb V^{m},
\\[2.5ex]
\disp
m((\vphi,\unl v),t^\veps\l)
=\min_{1\les j\les m}\{\veps+m(\vphi,\l),m(v_j,\l)\}, 
\quad \forall\,(\vphi,\unl v)\in\mbb W_{m}.
\end{array}
$$
The stability criterion for $\mbb W_m$ reads: a point 
$w=(\vphi,\unl v)$ is stable w.r.t. $(K^\times\times G,\,\th_m)$ 
if and only if 
\begin{equation}{\label{abc}}
\kern-.5ex\left\{
\begin{array}{lrlr}
(A)& 
m(\vphi,\l)\ges 0,\, m(\unl v,\l)\ges 0
&\Rightarrow&
\langle \chi_\ac,\l\rangle > 0;
\\[1ex] 
(B)&
1+m(\vphi,\l)\ges 0,\, m(\unl v,\l)\ges 0
&\Rightarrow& 
l+m\cdot \langle \chi_\ac,\l\rangle > 0;
\\[1ex]
(C)&
-1+m(\vphi,\l)\ges 0,\, m(\unl v,\l)\ges 0
&\Rightarrow& 
-l+m\cdot\langle \chi_\ac,\l\rangle > 0.
\end{array}
\right.
\end{equation}
Note that $\crl C_{(\vphi,\unl v)}
=\crl C_\vphi+\bigl(\mbb R\times\crl C_{\unl v}\bigr)$ 
for all $(\vphi,\unl v)\in\mbb W_m$; moreover, for $\unl v=(v_1,\ldots,v_m)$, 
then $\crl C_{\unl v}=\crl C_{v_1}+\ldots+\crl C_{v_m}$. 
We deduce that as {\em both} $m$ and $(\vphi,\unl v)\in\mbb W_m$ vary, 
there will be only {\em finitely many} dual cones: 
\begin{equation}{\label{fcs}}
\L_{(\vphi,\unl v)}
=
\L_\vphi\cap\bigl(\mbb R\times\L_{\unl v}\bigr)
=
\L_\vphi\cap\bigl(\mbb R\times(\L_{v_1}\cap\ldots\cap\L_{v_m})\bigr).
\end{equation}
We denote by $\L_1',...\,,\L_Z'$ the various intersections of 
$\L_1,...\,,\L_z$ defined above, corresponding to the 
{\em fixed} representation $G\rar\Gl(V)$.

\begin{proposition}{\label{prop:large-m}} 
Let us assume that the $G$-module $V$ has the property:\\ 
\centerline{
$(\bbV^m)^\sst(G,\chi_\ac)=(\bbV^m)^\s(G,\chi_\ac)\;$ for all $\;m\ges 1$.
} 
Then there is a constant $a_0(E)$ depending on $E$ 
such that for $\frac{m}{l}>a_0(E)$:\smallskip 

\centerline{$ 
\bigl(\mbb E\times\mbb V^{m}\bigr)^\s
(K^\times\times G,l\chi_t+m\chi_\ac)
=\bigl(\mbb E\sm\{0\}\bigr)\times \bigl(\mbb V^{m}\bigr)^\s(G,\chi_\ac). 
$}\smallskip

\nit Equivalently, $\chi_t+r\chi_\ac$ is an ample class on $\mbb P(\cal E)$ 
for $r>a_0(E)$. 
\end{proposition}

\begin{proof} 
`$\supset$'\quad 
Let $(\vphi,\unl v)\in (\mbb E\sm\{0\})\times (\mbb V^m)^\s(G,\chi_\ac)$. 
By definition, this means:$\;$
$
m(\unl v,\l)\ges 0 
\;\Rightarrow\;
\langle\chi_\ac,\l\rangle > 0.
$\\ 
The conditions $(A)$ and $(B)$ in \eqref{abc} are automatically fulfilled. 
We prove that for large $m$ the condition $(C)$ holds too. 
Let $\l_0$ be such that $m(\vphi,\l_0)\ges 1$ and $m(\unl v,\l_0)\ges 0$. 

Recall that only finitely many cones $\L_{\unl v}$ will appear when 
both $m$ and $\unl v\in(\mbb V^m)^\s$ vary. On each such cone, the linear 
function $\langle\chi_\ac,\,\cdot\,\rangle$ is strictly positive. 
We choose $a_1>0$ such that 
$
\langle\chi_\ac,\l\rangle\ges a_1|\l|,\; \forall\, 
\l\in\L_1'+\ldots+\L_Z'.
$

For fixed $\vphi$, the function $m(\vphi,\cdot)$ is piecewise linear. 
As $\vphi$ varies, $m(\vphi,\cdot)$ depends only on the weights of $T$ 
on $E$. Overall we find a constant $a_2(E)>0$ {\em independent of} 
$\vphi$ such that 
$|m(\vphi,\l)|\les a_2(E)\cdot|\l|$ for all $\l\in\cal X_*(T)$. 

Back to the proposition:
$$
\begin{array}{lr}
&a_2(E)\cdot|\l_0|\ges m(\vphi,\l_0)\ges 1
\;\Rightarrow\; 
|\l_0|\ges\frac{1}{a_2(E)}.
\\[2ex] 
\text{Hence:}\, 
&
-l+m\cdot \langle \chi_\ac,\l_0\rangle 
\ges 
-l+m\cdot a_1|\l_0|
\ges
-l+m\cdot \frac{a_1}{a_2(E)}.
\end{array}
$$
We conclude that for $\frac{m}{l}>\frac{a_2(E)}{a_1}$ the condition $(C)$ 
is satisfied. 

`$\subset$'\quad 
We prove that 
$$
\bigl(\mbb E\times\mbb V^{m}\bigr)^\us
(K^\times\times G,l\chi_t+m\chi_\ac)
\supset
\bigl(\mbb E\sm\{0\}\bigr)\times \bigl(\mbb V^{m}\bigr)^\us(G,\chi_\ac)
\;\text{ for }m\gg0.
$$ 
The conclusion follows from 
the hypothesis 
$(\mbb V^{m})^\sst(G,\chi_\ac)
=(\mbb V^{m})^\s(G,\chi_\ac)$. 

Recall from (\ref{fct-m}) that 
$\unl v\!\in\!(\mbb V^{m})^\us(G,\chi_\ac)$ if and only if 
$m(\unl v)\!<\!0$. The value $m(\unl v)$ is reached at the `worst' 
destabilizing $\l\in\cal X_*(G)$ (see \cite{ke2}). 
For variable $m$, there are only finitely many combinatorial strata 
in $(\mbb V^{m})^\us(G,\chi_\ac)$ ({\it c.f.} \eqref{fcs}), hence 
only finitely many possible values for $m(\unl v)$. It follows that  
$$
-\mu:=\max\bigl\{
m(\unl v)\mid m\ges 1,\,\unl v\in(\mbb V^{m})^\us(G,\chi_\ac)
\bigr\}<0.
$$
Now consider $(\vphi,\unl v)\in 
\bigl(\mbb E\sm\{0\}\bigr)\times \bigl(\mbb V^{m}\bigr)^\us(G,\chi_\ac)$, 
and its worst destabilizing $\l\in\cal X_*(G)$. After possibly moving 
$\unl v$ by an element in $G$, we may assume that $\l\in\cal X_*(T)$. 
Then holds: $m(\unl v,\l)\ges 0$, 
and $\frac{\langle\chi_\ac,\l\rangle}{|\l|}=m(\unl v)\les-\mu$. 

We distinguish the following cases: 

\nit-- If $m(\vphi,\l)=0\text{ resp.}>0$, then  
\eqref{abc}(A) and (B) imply that $(\vphi,\unl v)$ is $l\chi_t+m\chi_\ac$ 
unstable. 

\nit-- If $m(\vphi,\l)<0$, we normalize $\l$ such that $m(\vphi,\l)=-1$. 
We claim that $l+m\langle\chi_\ac,\l\rangle\les 0$ for $m$ large enough. 
Otherwise we deduce: 
$$
\phantom{xxxxxxxxxx}
\begin{array}{l}
\mu|\l|\les|\langle\chi_\ac,\l\rangle|<l/m
\\[1ex] 
1=|m(\vphi,\l)|\les a_2(E)\cdot|\l|
\end{array}
\biggr\}\,\Rightarrow\,
\frac{m}{l}<\frac{a_2(E)}{\mu}.
\phantom{xxxxxxxxxx}
$$
\end{proof}

%%%%%%%%%%%%%%%%%%%%%%%%%%%%%%%%%%%%%%%%%%%%%%%%%%%%%%%%%%%%%%%%%%%%%

\section{The nef vector bundles}
{\label{sct:nef-vb}}\setcounter{equation}{0}

In this section we define a finite set of `extremal' nef vector bundles, 
which will be the building blocks of the exceptional sequences. 
We continue the notations of the previous section. 
Consider the following Weyl group invariant cone:
\begin{equation}{\label{nef-cone}}
\crl N=\crl N(G,V):=\crl C_1\cap\ldots\cap\crl C_z
={(\L_1+\ldots+\L_z)}^\vee.
\end{equation}
When $G$ is a torus, $\crl N$ is the nef cone of the quotient, 
which is a toric variety. 
In our context, $\crl N$ can be viewed as the nef cone of 
$\mbb V^\sst(T,\chi_\ac)/T$. Its importance relies on the following:

\begin{proposition}{\label{prop:nef-vb}}
We consider a $G$-module $V$ which has the following property: 
$\bbV^\sst(G,\chi_\ac)=\bbV^\s(G,\chi_\ac)$. 
Let $E$ be a $G$-module, and $\cal E\rar Y$ be its associated 
vector bundle. 

Then $\cal E$ is nef if and only if 
all the weights of $T\!$ on $E\!$ belong to the cone $\crl N\!$. 
We call a module with this property {\em a nef module}.
\end{proposition}

\begin{proof}
\nit($\Leftarrow$) 
Let us assume that the weights of $E$ belong to $\crl N$. We prove that, 
on $\mbb P(\cal E^\vee)$, the class $\chi_t$ is nef, 
it means $\chi_t+r\chi_\ac$ is ample for all $r>0$. 
This translates into the following condition:  
$$
\bigl(\mbb E^\vee\times\mbb V\bigr)^\s
(K^\times\times G,\chi_t+r\chi_\ac)
=\bigl(\mbb E^\vee\sm\{0\}\bigr)\times
\mbb V^\s(G,\chi_\ac),
\quad\forall r>0.
$$
`$\supset$' The conditions \eqref{abc}(A) and (B) 
are trivially satisfied. We show that the case \eqref{abc}(C) 
does not occur. 

Take $(\psi,v)\in(\mbb E^\vee\sm\{0\})\times\mbb V^\s(G,\chi_\ac)$, 
and suppose that there is $\l_0$ with 
$m(\psi,\l_0)\!\ges 1$ and $m(v,\l_0)\!\ges 0$. Then  
$
\l_0\!\in{\rm int.}(\crl C_\psi)^\vee
\!\subset\!
-{\rm int.}\crl N^\vee
$
and also $\l_0\!\in\crl C_v^\vee\!\subset\!\crl N^\vee$. 
Contradiction. 

\nit`$\subset$' For shorthand, we denote $\cal S_L$ resp. $\cal S_R$ 
the left- and the right-hand-side above. 
Note that the quotient $\cal S_R\bigr/K^\times\kern-.7ex\times G$ 
exists, and equals $\mbb P(\cal E^\vee)$; 
let $Z:=\cal S_L/(K^\times\!\times G)$ be the quotient. 

By previous step, there is a morphism $\phi:\mbb P(\cal E^\vee)\rar Z$. 
Since $\phi$ is open and $\mbb P(\cal E^\vee)$ is projective, 
$\phi$ is surjective. Recall from \cite[Theorem 1.10]{mfk}, that 
$K^\times\!\times G$ acts with closed orbits on $\cal S_L$, and the 
quotient $\cal S_L\rar Z$ is geometric. 
Since $\mbb P(\cal E^\vee)\rar Z$ is surjective, the inclusion 
$\cal S_L\supset\cal S_R$ must be an equality. 
Otherwise we find closed orbits in $\cal S_L$, 
which are not contained in $\cal S_R$.

\nit($\Rightarrow$) Assume that $\cal E\rar Y$ is nef, that means 
$\chi_t$ is a nef class on $\mbb P(\cal E^\vee)$. By inspecting 
the conditions \eqref{abc} we deduce:
$$
\raise.15ex\hbox{$\not$}\exists\,\psi\in\mbb E^\vee\sm\{0\},\; 
v\in\bbV^\s(G,\chi_\ac),\;
\l\in\cal X_*(T)\text{ s.t. }
\left\{\begin{array}{l}
m(\psi,\l)\ges 1,\\ 
m(v,\l)\ges 0.
\end{array}\right.
$$
We choose $\psi=\vphi^\vee$, with $\vphi\in\mbb E$ of weight $\veps$. 
The previous condition implies: 
$\;\raise.15ex\hbox{$\not$}\exists\,
\l\in\cal X_*(T)$ such that $\langle\veps,\l\rangle<0$, and  
$\l\in\bigl(-\mbb R_+\veps+\crl N\bigr)^\vee.$
This happens only for $\veps\in\crl N$.
\end{proof}

There is also an effective procedure to produce `the smallest' 
such modules. Let us consider the set of weights:  
\begin{equation}{\label{n1}}
\crl N_1=\crl N_1(G,V):=
\left\{
\xi\,\biggl|\,
\begin{array}{l}
\mbb R_+\xi\text{ is an extremal ray of }\crl N,\\ 
\xi\text{ generates }\mbb R_+\xi\cap\cal X^*(T)
\text{ over }\mbb Z_{\ges0}
\end{array}
\right\}.
\end{equation}
It is a Weyl-invariant set, and therefore it makes sense considering the 
irreducible $G$-modules whose dominant weights belong to $\crl N_1$. 
These modules will be the building blocks for constructing exceptional 
sequences. We denote 
\begin{equation}{\label{eqn:nef-vb}}
\VB^+(Y):=
\left\{
E
\;\biggl|\;
\begin{array}{l}
\text{the dominant weight of the $G$-module $E$}
\\ 
\text{belongs to }\crl N_1
\end{array}\biggr.
\right\}.
\end{equation}
Equivalently, denote $\mfrak W_G^+$ the closure of the positive 
Weyl chamber of $G$. Then $\VB^+(Y)$ can be identified with 
$$\crl N_1^+(G,V):=\mfrak W_G^+\cap\crl N_1(G,V).$$

\begin{lemma}
The set $\VB^+(Y)$ is finite. For any $E\in\VB^+(Y)$, the weights of 
$T$ on $E$ belong to the cone $\crl N$.
\end{lemma}

\begin{proof}
As $\crl N_1$ is finite, $\VB^+(Y)$ is the same. Let $\xi$ be the 
dominant weight of $E$. The weights of $T$ on $E$ belong to the 
convex hull of the images of $\xi$ under the Weyl group. 
But all of them generate rays of $\crl N\!$. Hence 
the convex hull of the images of $\xi$ is contained in $\crl N\!$.
\end{proof}

\begin{proposition}{\label{prop:+comb}}
Let $M$ be an irreducible, nef $G$-module. 
Then there are $E_1,\ldots,E_n\in\VB^+(Y)$, and 
$c_1,\ldots,c_n\ges 1$ such that 
$M\!\subset\!
\mbox{$\overset{n}{\underset{j=1}\bigotimes}$}\Sym^{c_j}E_j.$ 
We say that $M$ is {\em a positive combination} of extremal 
nef modules.
\end{proposition}

\begin{proof}
Since the $G$-module $M$ is nef, its highest weight $\xi_M$ 
belongs to the cone $\crl N$. Then $\xi_M$ is a positive 
combination of $\xi_1,\ldots,\xi_n\in\crl N_1\,$:\\
\centerline{ 
$\xi_M=\mbox{$\overset{n}{\underset{j=1}\sum}$}c_j\xi_j,\;c_j\ges 1.$
}\\ 
Each $\xi_j$ is conjugated to some $\xi^+_j\in\crl N_1^+$, since 
the Weyl group acts transitively on the Weyl chambers. The irreducible 
$G$-module $E_j$ with highest weight $\xi^+_j$ belongs to $\VB^+(Y)$. 
Now observe that $\xi_M$ appears among the weights of 
$\overset{n}{\underset{j=1}\bigotimes}\Sym^{c_j}E_j$. 
Hence the whole module $M$ is contained in it.
\end{proof}

\begin{lemma}{\label{lm:chi}}
Consider the set $\VB^+(Y)$ of extremal nef vector bundles on $Y$, 
defined in \eqref{eqn:nef-vb}. Then the anti-canonical character 
$\chi_\ac(G,V)$ is a positive linear combination 
of $\det E$, with $E\in\VB^+(Y)$: 
$$
\chi_\ac=\mbox{$\underset{E\in\VB^+(Y)}\sum$}\kern-1ex 
m_{E}\!\cdot\det(E),\quad\text{with }m_E\ges 0.
$$
\end{lemma}

\begin{proof}
Let $\{\xi_j\}_j$ be the elements of $\crl N_1$. 
Since $\chi_\ac$ belongs to the interior of $\crl N$, there are 
positive numbers $c_j$ such that 
$
\chi_\ac\!=\!\underset{j}\sum c_j\xi_j
\!=\!\underset{j}\sum c_j\xi_j^\circ
+\underset{j}\sum c_j\xi_j'.
$
We decompose $\cal X^*(T)_{\mbb Q}
\!=\!\cal X^*(Z(G)^\circ)_{\mbb Q}\oplus\cal X^*(T')_{\mbb Q}$. 
Accordingly, each $\xi_j$ decomposes into 
$\xi_j=\xi_j^\circ+\xi_j'$, 
and each $\xi_j$ is conjugated to some $\xi_j^+\in\crl N_1^+$. 
Let $E_j\in\VB^+(Y)$ be the irreducible $G$-module with 
highest weight $\xi^+_j$. Note that $Z(G)^\circ$ acts on 
$E_j$ by the character $\xi_j^\circ$. 
Since $\chi_\ac$ is trivial on the semi-simple part of $G$, 
we deduce that $
\chi_\ac=\mbox{$\underset{j}\sum$} c_j\xi_j^\circ
=\mbox{$\underset{j}\sum$} 
\frac{c_j}{\dim E_j}\det E_j.$
\end{proof}

%%%%%%%%%%%%%%%%%%%%%%%%%%%%%%%%%%%%%%%%%%%%%%%%%%%%%%%%%%%%%%%%%%%%%

\section{Cohomological properties of nef vector bundles}
{\label{cohom-nef}}\setcounter{equation}{0}

In section \ref{sct:nef-vb} we have introduced the set of 
nef vector bundles associated to representations of $G$. 
In this section we are going to study their cohomological properties. 

\begin{theorem}{\label{thm:hq-nef}}
Let $E$ be a nef $G$-module. Then $H^q(Y,\cal E)=0$ for all $q>0$. 
\end{theorem}

\begin{proof}
Using the projection formula, 
$H^q(Y,\cal E)=H^q(\mbb P(\cal E^\vee),\cal O_{\mbb P}(1))$, and 
$\cal O_{\mbb P}(1)\rar \mbb P(\cal E^\vee)$ is a nef line bundle. 
The vanishing of the latter cohomology group is a consequence of 
the Hochster-Roberts theorem (see \cite{ke}). 
\end{proof}

We place ourselves in the following framework:
\begin{equation}{\label{rel-situation}}
\left\{\text{
\begin{minipage}{27.5em}
(i) There is a quotient group $H$ of $G$ with kernel $G_0$ 
(note that $G_0$ and $H$ are still reductive), 
and a quotient $H$-module $W$ of $V$ with kernel $V_0$, such that the 
natural projection 
$\pr^{\bbV}_{\mbb W}\!:\!\bbV\!\rar\!\mbb W$ has the property 
$\pr^{\bbV}_{\mbb W}\bigl(\;\bbV^\sst\bigl(G,\chi_\ac(G,V)\bigr)\;\bigr)
\subseteq\mbb W^\sst\bigl(H,{\chi_\ac}(H,W)\bigr).$\break  
We denote $Y\srel{\phi}{\rar}X$ the induced morphism. 
\\[1.5ex]
(ii) Both unstable loci have codimension at least two. 
\\[1.5ex]
(iii) $G$ and $H$ act freely on $\bbV^\sst\bigl(G,\chi_\ac(G,V)\bigr)$ 
and \phantom{MMM} 
$\mbb W^\sst\bigl(H,{\chi_\ac}(H,W)\bigr)$ respectively. 
\end{minipage}
}\right.
\end{equation}
Now let us study the positivity properties of direct images of nef 
vector bundles. 

\begin{lemma}{\label{lm:phiE}}
Suppose that we are in the situation \eqref{rel-situation}, and that 
$E$ is a $G$-module such that its associated vector bundle $\cal E\rar Y$ 
is nef. Then $\phi_*\cal E\rar X$ is a vector bundle, and it is 
associated to the $H$-module 
$\phi_*E:=\Mor(\bbV_0,E)^{G_0}=H^0(\bbV_0\invq_{\chi_\ac}G_0,\cal E)$. 
\end{lemma}

\begin{proof}
The restriction of $\cal E$ to the fibres of $\phi$ are nef. 
By applying theorem \ref{thm:hq-nef}, we obtain that $R^q\phi_*\cal E=0$ 
for all $q>0$, and therefore $\phi_*\cal E\rar X$ is locally free. 
Observe that both $V_0$ and $H$ are actually $G$-modules, and 
$V=V_0\oplus W$; the kernel $G_0$ is acting trivially on $W$. 
For an $H$-invariant open set $O\subset\mbb W$, holds:\smallskip 

\hspace{3.5em} 
$
\!\begin{array}[b]{ll}
H^0(O\invq H,\phi_*\cal E)
&=
H^0\bigl((\mbb V_0\times O)\invq G,\cal E\bigr)
=\Mor\bigl(\mbb V_0\times O,E\bigr)^G\!
\\
&=
\kern-.4ex{\bigl(\Mor(\mbb V_0\times O,E)^{G_0}\bigr)}^H
\kern-.8ex=\!
{\Mor\bigl(O,\Mor(\mbb V_0,E)^{G_0}\bigr)}^H
\kern-.5ex.
\end{array}
$
\end{proof}

\begin{theorem}{\label{thm:direct-image}}
Assume that {\rm($\!$}\ref{rel-situation}{\rm)} holds, and let $E$ be a 
nef $G$-module. Then the $H$-module $\phi_*E$ is still nef. 
(The direct image $\phi_*\cal E\rar X$ is a nef vector bundle.)
\end{theorem}  

Mourougane proves in \cite{mou} a similar statement for adjoint 
bundles. The proof below follows {\it ad litteram} his proof 
({\it loc.\,cit.} section 3), with the necessary changes. 

\begin{proof}
By lemma \ref{lm:phiE}, $\phi_*\cal E\rar X$ is locally free. 

\nit\unbar{Step 1}: 
Construct the tensor powers $(\phi_*\cal E)^{\otimes n}$.\\ 
Let $Y^{(n)}=Y\times_X\ldots\times_XY$ be the fibre product, and 
$\phi^{(n)}:Y^{(n)}\rar X$ be the projection. 
Note that the vector bundle 
$\cal E^{(n)}:=\cal E\times_X\ldots\times_X\cal E$ on $Y^{(n)}$ 
is nef. Its direct image is 
$\phi^{(n)}_*\cal E^{(n)}=(\phi_*\cal E)^{\otimes n}$. 
Moreover, $Y^{(n)}$ is the quotient of the affine space 
$\mbb V^{(n)}$ by the action of the group $G^{(n)}$, and 
$\cal E^{(n)}$ is associated to the $G^{(n)}$-module $E^{\oplus n}$: 

-- $V^{(n)}:=
\{(v_1,\ldots,v_n)\in V^{\oplus n}\mid \pr^V_W(v_1)=\ldots=\pr^V_W(v_n)\};
$

-- The group $G^{(n)}:=G\times_H\ldots\times_HG$ is still reductive.

\nit\unbar{Step 2}: 
Let $A\rar X$ be a very ample line bundle, associated to 
some character of $H$. Then 
$(\phi_*\cal E)^{\otimes n}\otimes A^{\dim X+1}$ is globally generated.\\ 
We replace $Y$ by $Y':=Y^{(n)}$, 
$\phi$ by $\phi':=\phi^{(n)}$, and $\cal E$ by $\cal E':=\cal E^{(n)}$. 

By the Castelnuovo-Mumford criterion, in order to prove that 
$\phi'_*\cal F\otimes A^{\dim X+1}$ is globally generated, 
it is enough to check that 
$H^q(X,\phi'_*\cal E'\otimes A^{\dim X+1-q})=0$ for all $q>0$. 
Since the higher direct images of $\cal E'$ vanish, the projection formula 
gives:\smallskip 
  
\centerline{$
H^q(X,\phi'_*\cal E'\otimes A^{\dim X+1-q})=
H^q(Y',\cal E'\otimes {(\phi')}^*A^{\dim X+1-q}). 
$}\smallskip

\nit But $Y'$ is still a quotient of an affine space, 
$\cal E'$ is associated to a nef $G$-module, and ${(\phi')}^*A$ 
corresponds to a nef character of $G$. We apply theorem 
\ref{thm:hq-nef} to $\cal E'\otimes {(\phi')}^*A^{\dim X+1-q}$, 
and deduce that its higher cohomology groups vanish. 

\nit\unbar{Step 3}: According to the previous step 
$(\phi_*\cal E)^{\otimes n}\otimes A^{\dim X+1}$ is globally 
generated for all $n>0$, and therefore $\phi_*\cal E$ is nef. 
\end{proof}

We use this result to describe more precisely the 
nef cone $\crl N(G,V)$. We consider the projective variety 
$$
\Flag(Y):=\O_G/B=\bigl(\O_G\times(G/B)\bigr)\bigr/G,
$$ and denote 
$\pi:\Flag(Y)\rar Y$ the projection. It is a $G/B$-fibre bundle over 
$Y\!$, justifying the notation $\Flag(Y)$. For any 
$\xi\in\cal X^*(T)=\cal X^*(B)$, we denote by 
$\cal L_\xi\rar\Flag(Y)$ the line bundle $(\O_G\times K)/B$, where 
$B$ acts on $K$ by $\xi$. 

\begin{corollary}{\label{LE}}
Let $\xi\in\cal X^*(T)$ be a dominant character, and let $E_\xi$ be 
the corresponding irreducible $G$-module. Then holds:
\begin{enumerate}
\item $\cal E_\xi=\pi_*\cal L_\xi$;

\item $\cal E_\xi\rar Y$ is nef if and only if 
$\cal L_\xi\in\Pic^+\bigl(\Flag(Y)\bigr):=$ the nef cone of $\Flag(Y)$.
\end{enumerate}
\end{corollary}

\begin{proof}
(i) The equality is a direct consequence of the Borel-Weil theorem, 
which says that $H^0(G/B,\cal L_\xi)=E_\xi$. 

\nit(ii) Assume that $\cal L_\xi$ is nef. The Borel-Weil theorem 
implies that the higher direct images $R^{>0}\pi_*\cal L_\xi=0$. 
By the same argument of the theorem \ref{thm:direct-image}, we deduce 
that $\cal E_\xi=\pi_*\cal L_\xi\rar Y$ is still nef. 

Conversely, assume that $E_\xi$ is nef, hence 
$\bbV^\sst(T,\chi_\ac)\subset\bbV^\sst(T,\xi)$. 
We claim that some tensor power of $\cal L_\xi$ is globally generated, 
and therefore $\cal L_\xi$ is nef. Let $B$ be the Borel subgroup of 
$G$ for which $\xi$ is dominant. Our hypothesis implies that 
$$
\begin{array}{r}
\bbV^\sst(G,\chi_\ac)=
\mbox{$\underset{g\in G}\bigcap$}g\bbV^\sst(T,\chi_\ac)
\subset
\mbox{$\underset{b\in B}\bigcap$}b\bbV^\sst(T,\chi_\ac)
\subset
\mbox{$\underset{b\in B}\bigcap$}
b\bbV^\sst(T,\xi)
\\=
\bbV^\sst(B,\xi).
\end{array}
$$
Observe that $B$ is solvable, {\em not reductive}, and therefore 
the standard invariant theory does not apply. 
The $B$-semi-stable locus $\bbV^\sst(B,\xi)$ is defined exactly as 
in \eqref{G-sst}, in terms of the algebra $K[\bbV]^{B,\xi}$. Its 
{\em finite generacy} has been proved by Grosshans (see {\it e.g.}
\cite[Corollary 9.5]{gross}). 

We deduce that for some $n>0$, $\bbV^\sst(B,\xi)$ 
can be covered by a finite number of sets $\{y\mid f(y)\neq 0\}$, 
with $f\in K[\bbV]^B_{\xi^n}$. Altogether, we find at each point 
$y\in\bbV^\sst(G,\chi_\ac)$ a function which is $(B,\xi^n)$-equivariant, 
and does not vanish at $y$. Hence $\cal L^n_{\xi}$ is globally 
generated. 
\end{proof}

\begin{corollary}{\label{cor:+comb}}
Suppose that \eqref{rel-situation} holds. 
Let $E$ be a nef $G$-module, and $M$ an irreducible $H$-submodule 
of $\phi_*E$. Then $M$ is a direct summand in a $H$-module 
of the form $\underset{F\in\VB^+(X)}\bigotimes\kern-1.2ex\Sym^{c_F}F.$
\end{corollary}

\begin{proof}
The push-forward $\phi_*\cal E\rar X$ is nef, and 
therefore all its weights belong to the cone $\crl N(H,W)$. 
We deduce that $M$ is nef too, and the conclusion follows from 
proposition \ref{prop:+comb}.
\end{proof}

\begin{example}
Consider the Grassmannian $X:=\Grs(K^m,d)$ of $d$-dimen\-sional quotients, 
and denote $\cal Q$ the tautological quotient on it. Note that the variety 
$\Flag(X)$ is the variety of full quotient flags of $\cal Q$. 
The cone $\mfrak W^+\cap\Pic^+\bigl(\Flag(X)\bigr)$ is generated by 
$d$ elements which correspond to the characters $\tau_1$, 
$\tau_1+\tau_2$,$\ldots$,$\tau_1+\ldots+\tau_d$ (here the $\tau_j$'s 
denote the obvious characters of the maximal torus in $\Gl(d)$). 

We deduce that for any nef $\Gl(d)$-module $F$, its associated vector 
bundle $\cal F\rar \Grs(K^m,d)$ is a direct summand in a tensor 
product of the form 

\centerline{
$\Sym^{c_1}(\cal Q)\otimes\Sym^{c_2}(\overset{2}\bigwedge\cal Q)
\otimes\ldots\otimes\Sym^{c_d}(\overset{d}\bigwedge\cal Q)$.}

\nit This is in agreement with the fact that this tensor product 
contains the Schur power $\mbb S^{\alpha}\cal Q$, where 
$\alpha=(\alpha_1\ges\ldots\ges\alpha_{d}\ges0)$, and the positive  
integers $c_j$ satisfy $\alpha_j=c_j+\ldots+c_d$ for $j=1,...,d$. 
\end{example}

%%%%%%%%%%%%%%%%%%%%%%%%%%%%%%%%%%%%%%%%%%%%%%%%%%%%%%%%%%%%%%%%%%%%%

\section{The main result: the absolute case}
{\label{sct:main}}\setcounter{equation}{0}

In this section we prove our first main result. 
We consider a $G$-module $V$, and the character $\chi_\ac=\chi_\ac(G,V)$. 
Assume that the codimension of the $\chi_\ac$-unstable locus is at 
least two, and $G$ acts freely on the semi-stable locus. 
It follows that 
$Y:=\bbV\invq_{\chi_\ac}G=\bbV^\sst(G,\chi_\ac)/G$ is a projective 
Fano variety. 
Observe that lemma \ref{lm:effective} implies that $\chi_\ac=\chi_\ac(G,V)$ 
is effective as soon as $m_\og > d_\og$ for all $\og\in\cal X$ 
(the result below does not require this hypothesis).

We define a Young diagram $\l$ of length $d$ to be an array 
of decreasing integers $(\l_1\ges\!\ldots\!\ges\l_d)$. 
We denote $\l_\mx:=\l_1$, $\l_\mn:=\l_d$, ${\rm length}(\l):=d$. 
For arrays consisting of positive integers, we visualize the Young 
diagrams, and the parameters as in the figure:\\ 
\centerline{\includegraphics{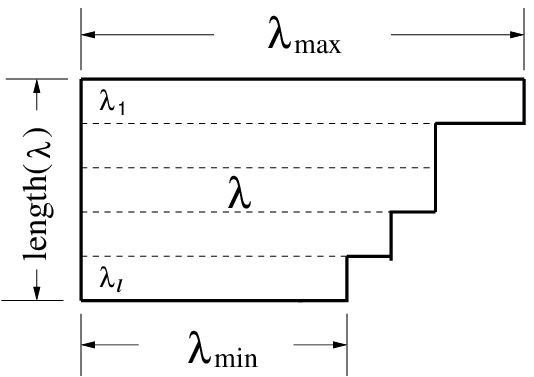}}
We introduce the following shorthand notation: 
for a Young diagram $\l$, let $\l\pm\fbox{$c$}$ 
be the diagram obtained by adding/subtracting 
the integer $c$ to/from the entries of $\l$. 
For a vector space $E$ and a Young diagram $\l$ of length $\dim E$, 
we will denote  $\mbb S^\l E$ its usual Schur power (for $\l_\mn\ges 0$), 
or $\mbb S^{\l- 
\raise.3ex\hbox{\tiny$\begin{array}{|c|}
\hline\kern-1.5ex \l_\mn\kern-1.5ex\\[.3ex] \hline
\end{array}$}
}\otimes(\det E)^{\l_\mn}$ (for arbitrary $\l$). 

For two positive numbers $m,d$ we define the following sets: 
\begin{eqnarray*}
\begin{array}{l}
\widetilde{\cal Y}_{d}:=
\text{ the set of Young diagrams $\l$ with }{\rm length}(\l)= d;
\\[2ex]
\cal Y_{m,d}:=
\bigl\{
\l\in\widetilde{\cal Y}_d\mid
0\les \l_\mn\les\l_\mx\les m
\bigr\};
\\[2ex] 
\cal Y_{d}:=\underset{m\ges 0}\bigcup\cal Y_{m,d}\,;
\quad
\cal Y^+_{d}:=
\underset{m\ges 0}\bigcup
\bigl\{
\l\in\cal Y_{m,d}\mid\l_d\ges{\rm length}
\bigl(
\l-\fbox{$\l_d$}
\,\bigr)\,
\bigr\}.
\end{array}
\end{eqnarray*}
Roughly speaking, our main result is that certain Schur powers of the 
extremal nef bundles on $Y$ form a strong exceptional sequence.

The main technical tool that will be used is the following cohomology 
vanishing theorem, proved by Manivel for K\"ahlerian varieties 
(see \cite{ma}), and Arapura for projective ones (see \cite{ar}).

\medskip\nit{\bf Theorem}
{\kern1ex\it Let $Y$ be a smooth projective variety, and 
$\{\cal E_1,\ldots,\cal E_N\}$ be a set of nef vector 
bundles over $Y$. Choose a set $\{\l^{(1)},\ldots,\l^{(N)}\}$ 
of Young diagrams such that 
$\l^{(j)}\in\cal Y_{\rk(\cal E_j)}^+$ for all $j$. 
Consider an ample line bundle $A\rar Y$. Then holds: 
$$
H^q\bigl(
Y,
\mbb S^{\l^{(1)}}\cal E_1\otimes\ldots\otimes\mbb S^{\l^{(N)}}\cal E_N
\otimes A\otimes\kappa_Y
\bigr)=0.
$$} 
\nit Next comes our first main result. 

\begin{theorem}{\label{thm:main}}
Let $V$ be a $G$-module such that $K[\mbb V]^T=K$. 
Assume that the unstable locus has codimension at least two, 
and that $G$ acts freely on $\bbV^\sst(G,\chi_\ac)$; 
we denote by $Y:=\bbV^\sst(G,\chi_\ac)/G$ the quotient. 
We consider the order $<_{\mfrak l}$ defined in \ref{defn:order2}. 

Let $E_1,\ldots,\kern-.1exE_N\!$ be the elements of 
$\VB^+\!(Y\kern-.1ex)$,$\kern-.2ex$ 
and denote $d_j\!:=\!\dim E_j$. 
We write $\chi_\ac=\overset{N}{\underset{j=1}\sum}m_j\cdot\det(E_j)$, 
with $m_j\ges 0$ as in lemma \ref{lm:chi}, and assume that all the 
numbers $m_j$ are integers. 
Consider the set 
$$
\begin{array}[t]{ll}
\euf{ES}(Y):=
&
\text{the set of all irreducible $G$-modules contained in}
\\[1ex] 
&
\mbb S^{\l^{\bullet}}E_\bullet
\!:=
\mbb S^{\l^{(1)}} E_1
\otimes\ldots\otimes
\mbb S^{\l^{(N)}} E_N, 
\text{ where }\l^{(j)}\!\in\!\cal Y_{m_j-d_j,d_j}.
\end{array}
$$
Then the vector bundles $\cal E\rar Y$ associated to the modules 
$E\in\euf{ES}(Y)$ form a strong exceptional sequence over $Y$ w.r.t. 
the order $<_{\mfrak l}$. 
\end{theorem}

\begin{proof}
The condition on $H^0(\Hom(\cal U',\cal U''))$ for two elements 
$\cal U',\cal U''\in\euf{ES}(Y)$ is implied by theorem \ref{thm:h000}. 
It remains to prove the vanishing of the higher cohomology groups. 
First of all we observe that, by definition, the vector bundles 
$\cal{U', U''}$ are direct summands of 
$\mbb S^{\l^\bullet}\cal E_\bullet$. 
Therefore it is enough to prove that vanishing of 
$H^q\bigl(Y,
\Hom(\mbb S^{\l^\bullet}\cal E,\mbb S^{\mu^\bullet}\cal E)
\bigr)$, $q>0$. 
Using the Littlewood-Richardson rules, we decompose  
$$
\Hom(\mbb S^{\l^\bullet}\cal E_\bullet,\mbb S^{\mu^\bullet}\cal E_\bullet)
=
\mbox{
${\underset{\alpha^\bullet=(\alpha^{(1)},\ldots,\alpha^{(N)})}\bigoplus}$
}
\kern-2ex\mbb S^{\alpha^\bullet}\cal E_\bullet\,,
$$
and observe that $\alpha^{(j)}\!=\!\bigl(m_j-d_j\ges
\alpha^{(j)}_1\ges\ldots\ges\alpha^{(j)}_{d_j}
\ges -m_j+d_j\bigr)$. For each direct summand holds: 
$$
\begin{array}{ll}
H^q\bigl(
Y,\mbb S^{\alpha^\bullet}\cal E_\bullet
\bigr)
&=
H^q\biggl(
Y,\kappa_Y\otimes\mbox{$\overset{N}{\underset{j=1}\bigotimes}$} 
\bigl(
\mbb S^{\alpha^{(j)}}\cal E_j\otimes \det(\cal E_j)^{m_j}
\bigr)
\biggr)
\\&
=
H^q\biggl(
Y,\kappa_Y\otimes\mbox{$\overset{N}{\underset{j=1}\bigotimes}$} \,
\mbb S^{\alpha^{(j)}+ 
\tiny\begin{array}{|c|}
\hline\kern-1.5ex m_j\kern-1.5ex\\ \hline
\end{array}
}
\cal E_j
\biggr).
\end{array}
$$
Note that 
$
\alpha^{(j)}+
\small\begin{array}{|c|}
\hline m_j\\ \hline
\end{array}
=
\underbrace{
\alpha^{(j)}+\small\begin{array}{|c|}
\hline -\alpha^{(j)}_{d_j}+d_j-1\\ \hline
\end{array}}_{:=\bar\alpha^{(j)}}
\,+\,
\small\begin{array}{|c|}
\hline \alpha^{(j)}_{d_j}+m_j-d_j+1\\ \hline
\end{array},$ 
and\vspace{-1.5ex}  
$$
\left\{\begin{array}{ll}
\bar\alpha^{(j)}_{d_j}=d_j-1\ges{\rm length}
\bigl(\bar\alpha^{(j)}-\small\fbox{$d_j-1$}
\,\bigr)
&\text{ and}
\\[1.5ex] 
\bar a_j:=\alpha^{(j)}_{d_j}+m_j-d_j+1\ges 1.
\end{array}\right.
$$
Since $E_1,\ldots,E_N$ are {\em all} the extremal nef bundles, 
it follows that the 
$A:=\overset{N}{\underset{j=1}\bigotimes}\det(\cal E_j)^{\bar a_j}$ 
is an ample line bundle over $Y$. The theorem cited above implies that 
the higher cohomology of $S^{\alpha^\bullet}\cal E_\bullet$ vanishes. 
\end{proof}

\begin{corollary}
Assume that the $G$-module $V$ has the property that the multiplicities 
$m_\og > d_\og$ for all $\og\in\cal X$. Then the exceptional 
sequence constructed above is formed by semi-stable vector bundles. 
\end{corollary}

\begin{proof}
It is an immediate consequence of the corollary \ref{thm:stab-bdl}.
\end{proof}

\begin{remark}{\label{rmk:length}}
It is important to observe that $\kappa_Y^{-1}$ is ample, and it 
becomes increasingly positive as we increase the multiplicities $m_\og $ 
of the isotypical components of $V$. 
It follows that the effect of increasing the $m_\og $'s is that 
of {\em simultaneously} increasing the dimension of the quotient, 
and that of the length of the exceptional sequence. 
In other words, for our construction we will always have a 
{\em lower bound} for 
$$\frac{\text{length of exceptional sequence on $Y$}}
{\text{Euler characteristic of } Y}.$$
Compare this construction with the one discussed in subsection 
\ref{ssect:AH}. 
\end{remark}

%%%%%%%%%%%%%%%%%%%%%%%%%%%%%%%%%%%%%%%%%%%%%%%%%%%%%%%%%%%%%%%%%%%%%

\section{The main result: the relative case}
{\label{sct:main2}}\setcounter{equation}{0}

Theorem \ref{thm:main} is too weak for fibred varieties. By applying 
it directly, one looses many terms of the exceptional sequences 
(see subsections \ref{ssect:kapranov} and \ref{ssect-A3}). 
The goal of this section is to address the relative case described 
in \eqref{rel-situation}. 
The additional hypothesis which will be imposed in 
\eqref{rel-situation2} may look overabundant, but in many 
concrete cases they are naturally fulfilled (especially 
for quiver representations). 

\begin{definition} 
Denote $T_0$ and $T_H$ the maximal tori of $G_0$ and $H$ respectively. 
The exact sequence 
$1\!\rar\! G_0\!\rar\! G\!\rar\! H\!\rar\!1$ induces a natural splitting 
$\cal X^*(T)_{\mbb Q}=\cal X^*(T_0)_{\mbb Q}\oplus\cal X^*(T_H)_{\mbb Q}$.

We will denote by $\crl N(G_0,V_0)$ respectively $\crl N(H,W)$ the nef 
cones of the $G_0$-module $V_0$ and $H$-module $W$, corresponding to  
$\chi_\ac(G_0,V_0)={\chi_\ac(G,V)|}_{G_0}$ and $\chi_\ac(H,W)$. 
\end{definition}

Throughout this section we will assume: 
\begin{equation}{\label{rel-situation2}}
\left\{\text{
\begin{minipage}{27.5em}
(i) The situation described in \eqref{rel-situation} holds.
\\[1ex]
(ii) $\crl N(G,V)=\crl N(G_0,V_0)+\crl N(H,W)$.\\ 
(We use the shorthand notation $\crl N=\crl N_0+\crl N_H$.)
\\[1ex]
(iii) The maximal torus $T_0\subset G_0$ has exactly $\dim T_0$ 
weights\\ on $V_0$.
\end{minipage}
}\right.
\end{equation}

\begin{remark}
Let us make a few comments related to the assumptions: 

\nit -- 
The condition (ii) means that there is a partition\smallskip 

\centerline{
$\VB^+(Y)=\VB^+(X)\,\dot\cup\,\VB^+(\text{fibre}).$}\smallskip 
 
\nit The set $\VB^+(X)$ can always be viewed as a subset of 
$\VB^+(Y)$ via the pull-back $\bbV\srel{\phi}{\rar}\mbb W$. 
What we assume is that the `extremal' nef bundles on 
the fibres extend to `extremal' nef bundles on the whole $Y$.
For shorthand, we will write $\VB^+_0:=\VB^+(\text{fibre})$.

\nit-- 
$T_0$ has always at least $\dim T_0$ linearly independent weights 
on $V_0$. The assumption (iii) is equivalent to any of the following: 

(iii$'$) For any $\xi\in\cal X^*(T_0)$, 
$\xi$ is $T_0$-nef on $V_0$ if and only if 
$\xi$ is $T_0$-effective on $V_0$;

(iii$''$) The quotient $\bbV_0\invq T_0$ is a product of projective 
spaces.
\end{remark}

\nit Observe that by lemma \ref{lm:chi}, we can express\smallskip 

\centerline{$
\begin{array}{l}
\chi_\ac(H,W)=\kern-1ex
\mbox{$\underset{F\in\VB^+(X)}\sum$}\kern-1.9ex
m_F\,\cdot\,\det F\;\;(m_F\ges 0),
\quad\text{and}
\\[2ex] 
\chi_\ac(G_0,V_0)=\kern-1ex
\mbox{$\underset{E\in\VB^+_0}\sum$}\kern-.4ex 
m_E\,\cdot\,\det E\;\;(m_E\ges 0).
\end{array}
$}

\begin{proposition}{\label{prop:tech}}
Assume that \eqref{rel-situation2} holds, and 
denote $d_F:=\dim F$, and $d_E:=\dim E$.

\nit{\rm(i)} Suppose that $(a_E)_{E\in\VB^+_0}$ and $(b_F)_{F\in\VB^+(X)}$ 
are integers having the following property: for all $q>0$, and all 
Young diagrams $\alpha^{E}\in\wtld{\cal Y}_{d_E}$ 
resp. $\beta^{F}\in\wtld{\cal Y}_{d_F}$, 
such that $\alpha^{E}_\mn\ges -a_E$ and $\beta^{F}_\mn\ges -b_F$, holds:
\begin{eqnarray}
{\label{eqn:a}}
H^q\biggl(
V_0\invq_{\chi_\ac(G_0,V_0)}G_0,
\mbox{$\underset{E\in\VB^+_0}\bigotimes$} 
\mbb S^{\alpha^{E}}\cal E
\biggr)=0,
\\ 
{\label{eqn:b}}
H^q\biggl(
X,\mbox{$\underset{F\in\VB^+(X)}\bigotimes$}\kern-1ex 
\mbb S^{\beta^{F}}\cal F
\biggr)=0.\hspace{5em}
\end{eqnarray}
Then $H^q\biggl(
Y,\underset{F\in\VB^+(X)}\bigotimes\kern-1.7ex
\phi^*\,\mbb S^{\beta^{F}}\cal F
\;\otimes
\underset{E\in\VB^+_0}\bigotimes\kern-.5ex
\mbb S^{\alpha^{E}}\cal E
\biggr)=0$ 
for all $q>0$, and for all Young diagrams 
$\beta^{F}\in\wtld{\cal Y}_{d_F}$ and $\alpha^{E}\in\wtld{\cal Y}_{d_E}$ 
with $\beta^{F}_\mn\ges -b_F$ and $\alpha^{E}_\mn\ges -a_E$ respectively.
\medskip

\nit{\rm(ii$_1$)} The condition \eqref{eqn:a} is fulfilled for 
$a_E:=m_E-d_E$, $\forall E\in\VB^+_0$. 

\nit{\rm(ii$_2$)} The condition \eqref{eqn:b} is fulfilled for 
$b_F:=m_F-d_F$, $\forall F\in\VB^+(X)$. 
\end{proposition}

\begin{proof}
$\!$(i) The hypothesis implies that the higher direct images of 
$\!\underset{E\in\VB^+_0}\bigotimes\kern-1ex\mbb S^{\alpha^{E}}\cal E$ 
vanish. By using the projection formula we deduce: 
\begin{eqnarray*}
H^q\biggl(
Y,\mbox{$\underset{F\in\VB^+(X)}\bigotimes$}\kern-1.7ex
\phi^*\mbb S^{\beta^{F}}\cal F
\;\otimes
\mbox{$\underset{E\in\VB^+_0}\bigotimes$}\kern-.5ex
\mbb S^{\alpha^{E}}\cal E
\biggr)\hspace{10em}
\\ 
=H^q\biggl(
X,\mbox{$\underset{F\in\VB^+(X)}\bigotimes$}\kern-1.7ex
\mbb S^{\beta^{F}}\cal F
\otimes\;
\phi_*\biggl(
\mbox{$\underset{\;E\in\VB^+_0}\bigotimes$}\kern-.5ex
\mbb S^{\alpha^{E}}\cal E
\biggr)\biggr).
\end{eqnarray*}
Let us write $\cal V^0:=
\underset{E\in\VB^+_0}\bigotimes\kern-.5ex
\mbb S^{\alpha^{E}}\cal E$, and decompose it into the direct sum 
corresponding to the irreducible $G$-modules appearing in the tensor 
product: $\cal V^0=\bigoplus\cal V^0_j$. The cohomology group breaks 
up into the direct sum of the `smaller' cohomology groups. 
For each component $\cal V^0_j$ there are two possibilities: 

\nit\unbar{Case 1} 
There is a weight of $T_0$ on $V^0_j$ which is not effective. 
In this case $\phi_*\cal V^0_j\!=\!0$ ({\it c.f.} theorem \ref{not-eff}), 
and we discard it from the direct sum. 

\nit\unbar{Case 2} 
All the weights of $T_0$ on $V^0_j$ are effective. In this case 
the hypotheses \eqref{rel-situation2} (ii)+(iii) imply that the weights 
of $V^0_j$ are nef, and therefore $\cal V^0_j\rar Y$ is nef itself. 
Using theorem \ref{thm:direct-image} and proposition \ref{prop:+comb}, 
we deduce that $\phi_*\cal V^0_j\rar X$ is nef, and is actually contained 
in $\underset{F\in\VB^+(X)}\bigotimes\kern-1.5ex\Sym^{c_F}\cal F$, 
with $c_F\ges 0$. 

The Littlewood-Richardson rules imply that the tensor product\smallskip 

\centerline{
$\mbox{$\underset{F\in\VB^+(X)}\bigotimes$}\kern-1.5ex
\mbb S^{\beta^{F}}\cal F\,\otimes
\mbox{$\underset{F\in\VB^+(X)}\bigotimes$}\kern-1.5ex
\Sym^{c_F}\cal F$}\smallskip
 
\nit breaks up into the direct sum of 
$\underset{F\in\VB^+(X)}\bigotimes\kern-1.5ex
\mbb S^{\bar\beta^{F}}\cal F$, 
with $\bar\beta^F_\mn\ges\beta^F_\mn+c_F\ges b_F$. 
By the hypothesis, their higher cohomology vanishes. 

(ii$_1$) Note that ${\kappa_{Y/X}^{-1}\bigr|}_{\rm fibre}
\!=\!\underset{E\in\VB_0^+}\sum\! m_E\cdot\det E$. 
Consider Young diagrams $(\alpha^E)_{E\in\VB^+_0}$ with 
$\alpha^E_\mn\ges d_E-m_E$ for all $E$. It holds: 
$$
\begin{array}{ll}
\mbox{$\underset{E\in\VB^+_0}\bigotimes$}\kern-1ex 
\mbb S^{\alpha^{E}}\cal E\otimes\,\kappa_{Y/X}^{-1}
\biggr|_{\rm fibre}
&=
\mbox{$\underset{\;E\in\VB^+_0}\bigotimes$}\kern-1ex 
\bigl(
\mbb S^{\alpha^{E}}\cal E\otimes(\det\cal E)^{m_E}
\bigr)
\biggr|_{\rm fibre}
\\ 
&=
\mbox{$\underset{E\in\VB^+_0}\bigotimes$}\kern-1ex 
{\mbb S^{{\alpha^{E}+ 
\tiny\begin{array}{|c|}
\hline\kern-1.5ex m_E\kern-1.5ex\\ \hline
\end{array}
}}\,\cal E}
\biggr|_{\rm fibre}, 
\end{array}
$$
and 
${\alpha^{E}+ 
\small\begin{array}{|c|}
\hline m_E\\ \hline
\end{array}
}
=
\underbrace{
\alpha^{E}+\small\begin{array}{|c|}
\hline -\alpha^{E}_{\mn}+d_E-1\\ \hline
\end{array}}_{:=\bar\alpha^{E}}
\,+\,
\small\begin{array}{|c|}
\hline \alpha^{E}_{\mn}+m_E-d_E+1\\ \hline
\end{array}\,$
with 
\vspace{-1ex}$$
\left\{\begin{array}{l}
\bar\alpha^{E}_{\mn}=d_E-1\ges{\rm length}
\bigl(\bar\alpha^{E}-\small\fbox{$d_E-1$}
\,\bigr)_{\,\mbox{,}}
%&
\\[1.5ex] 
\bar a_E:=\alpha^{E}_{\mn}+m_E-d_E+1\ges 1.
%&
\end{array}\right.
$$
Manivel and Arapura's theorem implies that 
$R^q\phi_*(\mbb S^{\alpha^\bullet}\cal E_\bullet)=0$, for all $q>0$.

(ii$_2$) Consider Young diagrams $(\beta^F)_{F\in\VB^+(X)}$ with 
$\beta^F_\mn\ges d_F-m_F$ for all $F$. Then holds: 
$$
\mbox{$\underset{F\in\VB^+(X)}\bigotimes$}\kern-2ex 
\mbb S^{\beta^{F}}\cal F\otimes\,\kappa_X^{-1}
=\!\!
\mbox{$\underset{F\in\VB^+(X)}\bigotimes$}\kern-2ex 
\bigl(
\mbb S^{\beta^{F}}\cal F\otimes(\det\cal F)^{m_F}
\bigr)
=\!\!
\mbox{$\underset{F\in\VB^+(X)}\bigotimes\,$}\kern-2ex 
{\mbb S^{{\beta^{F}+ 
\tiny\begin{array}{|c|}
\hline\kern-1.5ex m_F\kern-1.5ex\\ \hline
\end{array}
}}\,\cal F}. 
$$
We deduce the vanishing of the higher cohomology as in (ii$_1$). 
\end{proof}

\begin{theorem}{\label{thm:main2}} 
Assume that the conditions \eqref{rel-situation2} are satisfied, and 
that there are integers $(b_F)_{F\in\VB^+(X)}$ which fulfill the property 
\eqref{eqn:b}. Then the elements of the set ${\euf{ES}}(Y)$ defined 
below form a strong exceptional sequence of vector bundles over $Y$: 
\begin{eqnarray*}
\begin{array}{ll}
{\euf{ES}}(Y):=
&
\text{all the direct summands, corresponding to irreducible}
\\ 
&\text{$G$-modules contained in }
\\[0.5ex] &
\phi^*\bigl(\mbb S^{\l^{\bullet}}\cal F_\bullet\bigr)
\otimes
\mbb S^{\nu^{\bullet}}\cal E_\bullet
:=
\phi^*\bigl(
\mbox{\kern-2ex$\underset{\tiny F\in\VB^+(X)}\bigotimes$}
\kern-1.7ex\mbb S^{\l^{F}}\cal F\,
\bigr)
\otimes\kern-.2ex
\mbox{$\underset{E\in\VB^+_0}\bigotimes$\,}
\kern-.5ex\mbb S^{\nu^{E}}\cal E, 
\end{array}
\end{eqnarray*}
with $\l^F\in\cal Y_{b_F,\,d_F},$ and $\nu^{E}\in\cal Y_{m_E-d_E,\,d_E}$. 

Moreover, it holds: 
$H^q\biggl(
Y,\underset{F\in\VB^+(X)}\bigotimes\kern-1.7ex
\phi^*\,\mbb S^{\beta^{F}}\cal F
\,\otimes
\underset{E\in\VB^+_0}\bigotimes\kern-.5ex
\mbb S^{\alpha^{E}}\cal E
\biggr)\!=\!0$ 
for all $q>0$, and all Young diagrams $\beta^{F}\in\wtld{\cal Y}_{d_F}$ 
and $\alpha^{E}\in\wtld{\cal Y}_{d_E}$, 
with $\beta^{F}_\mn\ges -b_F$ and $\alpha^{E}_\mn\ges -(m_E-d_E)$ 
respectively.
\end{theorem}

\begin{proof}
Let $\cal U'$ and $\cal U''$ be two elements of ${\euf{ES}}(Y)$. 
The condition on the $H^0(\Hom(\cal U',\cal U''))$ follows again from 
theorem \ref{thm:h000}. 

It remains to prove the vanishing of $H^q(\Hom(\cal U',\cal U''))$, 
for $q\ges 1$. 
By using the Littlewood-Richardson rules, we deduce that 
$\Hom(\cal U',\cal U'')$ is direct summand in 
$\mbox{$\underset{\alpha^\bullet,\beta^\bullet}\bigoplus$}
\phi^*\bigl(\mbb S^{\beta^{\bullet}}\cal F_\bullet\bigr)
\otimes\mbb S^{\alpha^{\bullet}}\cal E_\bullet,$ with 
\vspace{-2ex}$$
\left\{\begin{array}{l}
\hskip2.75em b_F\ges\beta^F_\mx\ges\beta^F_\mn\ges-b_F,
\\[1ex] 
m_E-d_E\ges\alpha^E_\mx\ges\alpha^E_\mn\ges-m_E+d_E.
\end{array}\right.
$$
The conclusion of the theorem follows from proposition 
\ref{prop:tech}(ii$_1$). 
\end{proof}

\nit An immediate consequence of the previous theorem is the following: 

\begin{corollary}{\label{cor:tower}}
Assume the following assumptions hold:\smallskip

\begin{enumerate}  
\item There is a sequence of quotients 
$G\!\rar\! G_1\!\rar\!\ldots\!\rar G_k\!\rar\! 1$, 
with $\Gamma_{j}\!:=\!\Ker(G_{j}\!\rar\! G_{j+1})$. 
\item $V=W_1\oplus\ldots\oplus W_k$, 
where $W_j$ is a $G_j$-module for all $j$. 
We define $V_j:=W_j\oplus\ldots W_k$ for all $j$.  
\item The projections $\pr_j:V_j\rar V_{j+1}$ satisfy the conditions 
\eqref{rel-situation2}. The induced morphisms are denoted by  
$$
\phi_j\!:\!\mbb V_{j}\invq_{\!\chi_\ac(G_{j},V_{j})} G_{j}
\rar \mbb V_{j+1}\invq_{\!\chi_\ac(G_{j+1},V_{j+1})} G_{j+1},\; 
\text{for all }\,1\les j\les k-1.
$$
\end{enumerate}\smallskip

\nit Let us write $\chi_\ac(\Gamma_j,W_j)
=\kern-.9ex
\underset{E\in\VB^+(\mbb W_j/\!/ \Gamma_j)}\sum\kern-3ex 
m_{j,E}\cdot\det E$ ({\it c.f.} \ref{lm:chi}), and denote 
$\VB^+_j:=\VB^+(\mbb W_j/\!/ \Gamma_j)$. 

Then the elements of the set ${\euf{ES}}(Y)$ defined below 
form a strong exceptional sequence of vector bundles over $\bbV\invq H$: 
$$
\begin{array}{ll}
{\euf{ES}}(Y):=
&
\text{all the direct summands, corresponding to irreducible}
\\ 
&
\text{$G$-modules contained in }
\ouset{j=1}{k}{\bigotimes}
\biggl(
\mbox{$\underset{\tiny E\in\VB^+_j}\bigotimes$}
\mbb S^{\alpha^{j,E}}\cal E
\biggr),
\\[1ex] &
\text{with }\;\alpha^{j,E}\in\cal Y_{m_{j,E}-d_E\,,\,d_E}.
\end{array}
$$
Assume moreover that the multiplicity condition in corollary 
\ref{cor:stab-bdl} is fulfilled. Then $\euf{ES}(Y)$ consists of 
semi-stable vector bundles over $Y$.
\end{corollary}

%%%%%%%%%%%%%%%%%%%%%%%%%%%%%%%%%%%%%%%%%%%%%%%%%%%%%%%%%%%%%%%%%%%%%

\section{Examples}
{\label{sct:expl}}\setcounter{equation}{0}

In this section we are going to present a few particular cases, 
in order to illustrate the general discussion. 
We concentrate on quiver varieties because they are a source of 
infinitely many examples, and are also very convenient: 
for generic choices of the dimension vector, the semi-stability and 
stability concepts agree. Therefore the quotients which will appear 
are geometric, as we wish. Even more remarkably, the procedure of 
constructing exceptional sequences of vector bundles over quiver 
varieties is {\em almost algorithmic}. 

Let $Q=(Q_0,Q_1,h,t)$ be a quiver, and $\unl d=(d_q)_{q\in Q_0}$ 
be a dimension vector. We adopt the following convention: 
suppose that $q,q'$ are two vertices, and there is (at least) one 
arrow from $q$ to $q'$; then we draw {\em only one} arrow $a$, and 
we denote by $m_a$ its {\em multiplicity} (that is how many times 
the arrow is repeated). In other words, we consider the group 
$G=\kern-.5ex
\underset{q\in Q_0}\btimes\kern-.5ex\Gl(d_q)$, and the $G$-module 
$V=\underset{a\in Q_1}\bigoplus\kern-.5ex
{\Hom(K^{d_{t(a)}},K^{d_{h(a)}})}^{\oplus m_a}\!$. 
The construction of exceptional sequences involves the following steps:
\smallskip

\nit\unbar{Step 1} Compute the anti-canonical character: 
$$
\begin{array}{rl}
\chi_\ac
=&
\mbox{$\underset{a\in Q_1}\sum$} m_a\cdot
\bigl(
d_{t(a)}\det_{h(a)}-d_{h(a)}\det_{t(a)}
\bigr)
\\[2ex]
=&
\mbox{$\underset{q\in Q_0}\sum$} 
\biggl(
\mbox{$\underset{a\in Q_1^{\rm in}(q)}\sum$}
\kern-1ex m_ad_{t(a)}
-
\mbox{$\underset{a\in Q_1^{\rm out}(q)}\sum$}
\kern-1ex m_ad_{h(a)}
\biggr)\cdot\det_q.
\end{array} 
$$
Note that the multiplicative group, embedded diagonally in $G$, 
acts trivially on $V$, and the quotient $G/{(K^\times)}_{\rm diag}$ 
acts effectively on $V$. Moreover, for generic choices of the 
multiplicities $m_a$ (w.r.t. the dimension vector $\unl d$),  
the $\chi_\ac$-semi-stable locus of $\bbV$ coincides with the 
stable locus (see {\it e.g.} \cite[proposition 3.1]{king}). 
For such a generic choice, there is a natural `Euler sequence' over 
the quotient $Y$:
$$
\vspace{-1.5ex}
0
\lar
\cal O_Y^{^{\oplus\dim\hat G}}
\lar
\mbox{$\underset{a\in Q_1}\bigoplus$}
{\cal Hom\bigl(\cal E_{t(a)},\cal E_{h(a)}\bigr)}^{\oplus m_a}
\lar
T_Y
\lar 0.
$$
It follows that the anti-canonical class of the quotient is 
$\kappa_Y^{-1}=\chi_\ac$.\smallskip

\nit\unbar{Step 2} 
It consists in determining the `extremal bundles' in the set $\VB^+(Y)$ 
(see \eqref{eqn:nef-vb}), and expressing $\chi_\ac$ as a positive 
combination of their determinants (see lemma \ref{lm:chi}). 
Actually this step is responsible for the use of the word 
`almost' above: the computation of the extremal nef bundles is 
algorithmic, but involves  the maximal torus of $G$, and is 
therefore tedious.\smallskip 

\nit\unbar{Step 3} Denote $\cal E_1,\ldots,\cal E_N$ the extremal 
bundles above, and take tensor products of their Schur powers 
$
\mbb S^{\l_1,\ldots,\l_N}\cal E\kern-.5ex
:=
\mbb S^{\l_1}\cal E_1\otimes\ldots\otimes\mbb S^{\l_N}\cal E_N.
$
The third step consists in determining the sizes of the Young 
diagrams $\l_1,\ldots,\l_N$ which fulfill the requirements of 
theorem \ref{thm:main}. 
\smallskip

\nit\unbar{Step 4 (Optional)} Search for fibrations coming from 
a sub-quiver. More precisely, we are looking 
for a sub-quiver $R\subset Q$ having the property: 
$$
\begin{array}{rcr}
\forall\,(A_a)_{a\in Q_1}\in\bbV^\sst\bigl(G,\chi_\ac(V)\bigr)
&\quad\Longrightarrow\quad&
(A_a)_{a\in R_1}\in\mbb W^\sst\bigl(H,\chi_\ac(W)\bigr),
\\[1.5ex] 
G=\underset{v\in Q_0}\prod\Gl(v)
&&
H=\underset{v\in R_0}\prod\Gl(v).
\end{array}
$$
Here $V$ and $W$ denote the representation spaces of $Q$ and $R$ 
respectively. 
In such a situation there is a natural projection map $Y\rar X$ 
between the corresponding quotients. 
Moreover, if $R$ is chosen appropriately, the numerous hypotheses 
in \eqref{rel-situation2} are naturally fulfilled. 

Very often one obtains better bounds for the sizes of the Young 
diagrams involved in the Schur powers  than those which are 
obtained by applying the step 3 directly 
(see subsections \ref{ssect:kapranov} and \ref{ssect-A3} below).

%%%%%%%%%%%%%%%%%%%%%%%%%%%%%%%%%%%%%%%%%%%%%%%%%%%%%%%%%%%%%%%%%%%%%%

\subsection{Kapranov's construction}{\label{ssect:kapranov}}
Let us start by reviewing Kapranov's examples of tilting bundles 
over the Grassmannian, and over the flag variety for $\Gl(m)$. 
We show that by using our approach we automatically recover the vector 
bundles which appear in the tilting objects constructed by Kapranov 
over the Grassmannian, and over partial flag manifolds. 

They are the quiver varieties associated respectively to:
$$
\entrymodifiers={++[o][F-]}
\xymatrix@+.7em@R=.9em{
m
\ar[r]^-{B}
&
*++[o][F=]{_{\phantom{i}}d_{\phantom{i}}} 
&*\txt{}&*\txt{}
&*\txt{\kern-1ex with $m>d$.}
\\
m
\ar[r]^-{A_{k}}
&
*++[o][F=]{d_k} 
\ar[r]^-{\!A_{{k-1}}} 
&*\txt{$\;$\ldots}\ar[r]^-{A_{1}}&
*++[o][F=]{d_1} 
&
*\txt{\kern-1ex with $m>d_k>\ldots>d_1$.}
}
$$
A doubled circle means that the corresponding linear group acts 
at that entry (we have factored out the diagonal $K^\times$-action).

%%%%%%%%%%%%%%%%%%%%%%%%%%%%%%%%%%%%%%%%%%%%%%%%%%%%%%%%%%%%%%%%%%%%%%

\subsubsection{The case of the Grassmannian} 
Let us consider the Grassmannian $Y:=\Grs(\mbb C^m,d)$ of $d$-dimensional 
quotients of $K^m$. Its anti-canonical class is 
$\kappa_{\Grs(K^m,d)}^{-1}=(\det\cal Q)^m$, where $\cal Q$ denotes the 
universal quotient bundle. The cone $\crl N$ is generated by the 
characters $t_1,\ldots,t_d$ of $\Gl(d)$, and $\crl N_1^+=\{t_1\}$. 
Hence the set $\VB^+(Y)$ of extremal nef bundles $\VB^+(Y)$ consists 
of $\cal Q$ only. 

Theorem \ref{thm:main} says that the elements of the set 
$\{\mbb S^\l\cal Q\mid \l\in\cal Y_{m-d,d}\}$ 
form a strong exceptional sequence of vector bundles on $\Grs(K^m\!,d)$. 
Indeed, this is what Kapranov proves in {\it loc.\,cit.}. Let us remark 
that he actually proves that they form a tilting sequence. 
\medskip

%%%%%%%%%%%%%%%%%%%%%%%%%%%%%%%%%%%%%%%%%%%%%%%%%%%%%%%%%%%%%%%%%%%%

\subsubsection{The case of flag manifolds}
We denote by $\mbb F_k:=\Flag(K^m,d_k,\ldots,d_1)$ 
the variety of quotient $k$-flags of $K^m$. Let 
$\cal Q_1,\ldots,\cal Q_k$ be the tautological quotient bundles over 
$\mbb F_k$ with $\rk\cal Q_j\!=\!d_j$. 

The anti-canonical class is $\kappa_{\mbb F_k}^{-1}
\!\!=\!\!
\hbox{$\overset{k}{\underset{j=1}\bigotimes}$}
{(\det\cal Q_j)}^{d_{j+1}-\,d_{j-1}}\!\!.$ 
The cone $\crl N$ is generated by the characters 
$t^{(j)}_1,...,t^{(j)}_{d_j}$, $j=1,...,k$, and 
$\crl N_1^+=\{t^{(1)}_1,\ldots,t^{(k)}_1\}$. 
We deduce that $\VB^+(\mbb F_k)=\{\cal Q_1,\ldots,\cal Q_k\}$. 
By applying theorem \ref{thm:main} directly, we obtain that the 
elements of 
$$
\left\{
\begin{array}{l}
\mbb S^{\l_\bullet}\cal Q_\bullet^\vee:=
\mbb S^{\l_k}\cal Q_k\otimes\ldots\otimes\mbb S^{\l_1}\cal Q_1,
\quad
\l_\bullet=(\l_k,\ldots,\l_1),
\\[1ex]
\text{with }\l_\bullet\in
\cal Y_{m-d_k-d_{k-1},d_k}
\times\ldots\times 
\cal Y_{d_3-d_2-d_{1},d_2}\times\cal Y_{d_2-d_1,d_1}\!
\end{array}
\right\}
$$
form a strong exceptional sequence over $\mbb F_k$. 
The problem is that these bounds are very weak, 
and this set can be empty! 

At this point Step 4 becomes useful. There is a natural 
projection from the $k$-flag onto the $(k-1)$-flag variety
$$
\mbb F_k\srel{\phi}{\lar}\mbb F_{k-1},\quad
[A_k,\ldots,A_2,A_1]\lmt[A_k,\ldots,A_2].
$$
One checks easily that all the conditions of \eqref{rel-situation2} 
are fulfilled. By applying corollary \ref{cor:tower} we deduce that 
the elements of the set 
$$
\left\{
\begin{array}{l}
\mbb S^{\l_\bullet}\cal Q_\bullet^\vee:=
\mbb S^{\l_k}\cal Q_k^\vee\otimes\ldots\otimes\mbb S^{\l_1}\cal Q_1^\vee
\\[1ex] 
\text{with }\l_\bullet=(\l_k,\ldots,\l_1)
\in
\cal Y_{m-d_k,d_k}\times\ldots\times\cal Y_{d_2-d_1,d_1}
\end{array}
\right\}
$$
form a strong exceptional sequence of vector bundles over $\mbb F_k$.

%%%%%%%%%%%%%%%%%%%%%%%%%%%%%%%%%%%%%%%%%%%%%%%%%%%%%%%%%%%%%%%%%

\subsection{$A_3$-type quiver with multiple arrows}
{\label{ssect-A3}} 
Interesting phenomena occur already for $A_3$-type quivers, as soon as 
we increase the multiplicities of the arrows. Consider the quiver 
$$\begin{array}{l}
\entrymodifiers={++[o][F-]}
\xymatrix@+.9em{
m
\ar[rr]^-{B}
&*\txt{}&
*++[o][F=]{d_2} 
\ar[rr]^-{{\mbb A}=(A_1,\ldots,A_{n})} 
&*\txt{}&
*++[o][F=]{d_1} 
&
*\txt{with $m>d_2>d_1$,}
}
\\[2ex] 
V={(K^{d_2})}^{\oplus m}\oplus{\Hom(K^{d_2},K^{d_1})}^{\oplus n},
\quad 
G=\Gl(d_1)\times\Gl(d_2).
\end{array}$$
Let $Y:=\bbV\invq_{\chi_\ac}G$ be the corresponding quiver variety. 
The flag variety $\Flag(K^m,d_2,d_1)$ corresponds to the case $n=1$. 
We denote the vector bundles over $Y$ associated to the $G$-modules 
$K^{d_1}$ and $K^{d_2}$ by $\cal E_1$ and $\cal E_2$ respectively. 
The anti-canonical character is 
$$
\chi_\ac=nd_2\cdot\det_1+(m-nd_1)\cdot\det_2
\!=\!
n\cdot\bigl[
d_2\cdot\det\cE_1+(r-d_1)\cdot\det\cE_2
\bigr],
\; r:=\!\frac{m}{n}.
$$
We are going to see that the effect of introducing the parameter $n$ 
is that of obtaining several types of quotients. 
Observe that for generic choices of $m$ and $n$, the semi-stable 
and the stable loci coincide; this happens for 
$$
\begin{array}{c}
\text{gcd}(nd_2,m-nd_1)=\text{gcd}(nd_2,m-nd_1,md_2)=1.
\end{array}
$$
For details about semi-stability criteria for quiver representations, 
the reader may consult \cite{king}.

\subsubsection{Case $r\!>\!d_1$} 
The $\chi_\ac$-semi-stability condition for 
$(B,\mbb A)\!\in\!\mbb V$ is: 
\begin{eqnarray}
\begin{array}{l}
\left\{
\begin{array}[c]{l}
U_2\subset K^{d_2}\text{ and }U_1\subset K^{d_1}\text{ s.t. }
\mbb A(U_2)\subset U_1
\\ 
\dim(U_2)=d_2'\text{ and }\dim(U_1)=d_1'
\end{array}
\right\}
\\[3ex] 
\hspace{3ex}\Longrightarrow
d_2d_1'+(r-d_1)d_2'\ges rd_2\text{ for }(d_2',d_1')\neq(d_2,d_1).
\end{array}
\end{eqnarray}
The set of extremal nef vector bundles is 
$\VB^+(Y)=\{\cal E_1,\cal E_2\}$, and the anti-canonical class is 
$\kappa_Y^{-1}=(\det\cal E_2)^{m-nd_1}\otimes(\det\cal E_1)^{nd_2}$. 
Theorem \ref{thm:main} implies that the elements of the set 
$$
\{
\mbb S^\l\cal E_1\otimes
\mbb S^\mu\cal E_2
\mid 
\l\in\cal Y_{nd_2-d_1,d_1}\text{ and }
\mu\in\cal Y_{m-nd_1-d_2,d_2}
\}
$$
form a strong exceptional sequences of vector bundles over $Y$. 

We illustrate again the role of Step 4 described at the 
beginning of this section: by using an appropriate fibre bundle 
structure on $Y$, we will increase the number of elements in 
the exceptional sequence. 

Observe that both $B$ and $\mbb A\in\Hom(K^{d_2}\otimes K^n,K^{d_1})$ 
are surjective, for any $\chi_\ac$-semi-stable point $(B,\mbb A)$. 
Indeed: by inserting $d_1'=d_1$ we obtain $d_2'\ges d_2$, 
and by inserting $d_2'=d_2$ we obtain $d_1'\ges d_1$. 
It follows that there is a natural projection $\phi:\!Y\!\rar\Grs(K^m,d_2)$, 
whose fibres are isomorphic to $\Grs(K^{nd_2},d_1)$. 
The group $\Gl(m)\times\Gl(n)$ acts on $Y$, and the projection is 
equivariant for the $\Gl(m)$-action. 
However $Y$ is not the $2$-flag variety. 

We observe that the projection $V\rar\Hom(K^m\!,K^{d_2})$ fulfills the 
conditions \eqref{rel-situation2}, and moreover 
$\VB^+\bigl(\Grs(K^m\!,d_2)\bigr)=\{\cal E_2\}$, and 
$\VB^+_0=\{\cal E_1\}$. Applying corollary \ref{cor:tower} to 
$\phi$ we deduce that the elements of the following set form a strong 
exceptional set of vector bundles over $Y$:
\begin{eqnarray}
\{
\mbb S^\l\cal E_1\otimes
\mbb S^\mu\cal E_2
\mid 
\l\in\cal Y_{nd_2-d_1,d_1}\text{ and }
\mu\in\cal Y_{m-d_2,d_2}
\}.
\end{eqnarray}

\subsubsection{Case $r\!<\!d_1$} 
The $\chi_\ac$-semi-stability condition for $(B,\mbb A)\!\in\!\mbb V$ is: 
\begin{eqnarray}{\label{ss2}}
\begin{array}{l}
\left\{\kern-.5ex
\begin{array}[c]{l}
U_2\subset K^{d_2}\text{ and }U_1\subset K^{d_1}\text{ s.t. }
\mbb A(U_2)\subset U_1
\\ 
\dim(U_2)=d_2'\text{ and }\dim(U_1)=d_1'
\end{array}\kern-.5ex
\right\}
\\[3ex] 
\hspace{3ex}
\Longrightarrow
\left\{\begin{array}{rll}
\rm (i)\;&
d_2d_2'-(d_1-r)d_1'\ges 0&\text{for }(d_2',d_1')\neq(0,0),
\\ 
\rm (ii)\;&
d_2d_1'-(d_1-r)d_2'\ges rd_2&\text{for }(d_2',d_1')\neq(d_2,d_1).
\end{array}\right.
\end{array}
\end{eqnarray}
Now we determine the extremal nef vector bundles. 
Since $r-d_1<0$, the situation differs from the previous case; 
now we will have $\VB^+(Y)=\bigl\{\cal E_2, \cal H\bigr\}$, 
with $\cal H:=\Hom(\cal E_2,\cal E_1)$. 
We express the anti-canonical class as a positive combination 
of the extremal bundles: 
$\kappa_Y^{-1}=(\det\cal E_2)^m\otimes(\det\cal H)^n$. 
Theorem \ref{thm:main} implies that 
\begin{eqnarray}
\{\mbb S^\l\cal E_2\otimes\mbb S^\mu\cal H
\mid
\l\in\cal Y_{m-d_2,d_2}\text{ and }
\mu\in\cal Y_{n-d_1d_2,d_1d_2}
\}
\end{eqnarray}
is a strong exceptional sequence of vector bundles over $Y$. 

Let us interpret the result. We consider the sub-quiver formed by the 
two rightmost vertices, 
and let $W:=\Hom(K^{d_2}\otimes K^n,K^{d_1})$ be its representation 
space. The anti-canonical character is $\chi_\ac(W)=d_2\det_1-d_1\det_2$. 
The symmetry group which is acting (effectively) is 
$G/{(K^\times)}_{\rm diag}$. The condition 
\eqref{ss2}(i) implies that if $(B,\mbb A)$ is $\chi_\ac$-semi-stable, 
then $\mbb A$ is $\chi_\ac(W)$-semi-stable. Hence there is a natural 
morphism 
$$Y\srel{\phi}{\lar} X:=
\Hom(K^{d_2}\otimes K^n\!,K^{d_1})\invq_{\chi_\ac(W)}\,G,
$$ 
which is a projective bundle, with fibre isomorphic to $\mbb P(K^{md_2})$.  
The conditions \eqref{rel-situation2} are fulfilled, and 
we may apply corollary \ref{cor:tower}. However, in this case we 
do not improve the previous bound.

%%%%%%%%%%%%%%%%%%%%%%%%%%%%%%%%%%%%%%%%%%%%%%%%%%%%%%%%%%%%%%%%%%%

\subsection{Altman and Hille's examples}{\label{ssect:AH}}
In the article \cite{ah} the authors present the following construction: 
consider a quiver $Q$ without oriented cycles, and a thin and faithful 
representation space $V$ of it. This means that the dimension vector 
of the representation space is $\unl d={(1)}_{q\in Q_0}$, and 
the symmetry group which is acting is the torus 
$T=\underset{q\in Q_0}\prod K^\times\bigl/{(K^\times)}_{\rm diag}$.
\medskip 

\nit{\bf Theorem} (Altmann, Hille) {\it 
Assume that $\bbV^\sst(T,\chi_\ac)=\bbV_{(0)}^\s(T,\chi_\ac)$. Then the 
tautological line bundles ${(\cal L_q)}_{q\in Q_0}$ form an exceptional 
sequence over the toric variety $Y:=\bbV^\sst(T,\chi_\ac)/T$. 
}\medskip

We wish to remark that this construction fits into a more 
general framework: we consider a quiver $Q=(Q_0,Q_1,h,t)$ without 
oriented cycles, and we fix a dimension vector $\unl d=(d_q)_{q\in Q_0}$; 
we denote $V$ the corresponding representation space. 
For $m\ges 1$, we denote $Q^{(m)}$ the quiver obtained from $Q$ by 
multiplying each arrow $m$ times. 
The representation space of $Q^{(m)}$ with dimension vector $\unl d$ 
is $V^m$, and the symmetry group which is acting 
is $G=\underset{q\in Q_0}\prod \Gl(d_q)\bigl/K^\times_{\rm diagonal}$. 

\begin{proposition}{\label{prop:general-AH}}
Assume that ${(\bbV^m)}^{\sst}(G,\chi_\ac)={(\bbV^m)}^{\s}(G,\chi_\ac)$, 
and denote $Y_m$ the quotient by the $G$-action. 
For $q\in Q_0$, we denote $\cal E_q$ the tautological bundle over $Y_m$, 
associated to $G\rar\Gl(d_q)$. 

Then there is a constant $m(Q)\ges 1$ such that for all $m> m(Q)$, 
the set $\{\cal E_q\}_{q\in Q_0}$ is a strong exceptional sequence of 
vector bundles over $Y_m$ (with respect to an appropriate ordering). 
Moreover, these vector bundles are semi-stable. 
\end{proposition}

\begin{proof}
For two vertices $p,q\in Q_0$, let 
$E_{pq}:=\Hom(E_p, E_q)$, and $\cal E_{pq}$ the associated vector 
bundle over $Y_m$, and let $e_{pq}:=\dim E_{pq}=d_pd_q$. 

The condition on 
$H^0(Y_m,\cal E_{pq})$ follows from theorem \ref{thm:h000}. 
It remains to prove the vanishing of the higher cohomology. 
We compute $H^{n-i}\bigl(Y_m,\cal E_{pq}\bigr)$, $n=\dim Y$, 
by using the relative duality for 
$\mbb P(\cal E_{qp})\srel{\pr}{\rar}Y_m$; it equals: 
$$
H^{(e_{pq}-1)+i}\bigl(\mbb P(\cE_{qp}),
\pr^*(\kappa_{Y_m}\otimes(\det\cE_{pq})^{-1})
\otimes\cO_{\mbb P(\cE_{qp})}(-e_{pq}-1)
\bigr)^\vee.
$$
The Kodaira vanishing theorem implies that 
$H^j(Y_m,\cal E_{pq})$ vanishes for all $j\ges 1$, 
as soon as 
$\pr^*(\kappa_{Y_m}^{-1}\otimes(\det\cE_{pq}))
\otimes\cO_{\mbb P(\cE_{qp})}(e_{pq}+1)$ 
is ample over $\mbb P(\cE_{qp})$. 
By proposition \ref{prop:large-m}, there is 
a number $m_{pq}$ such that this property holds for all $m>m_{pq}$. 
Consider now $m(Q):=\max\{m_{qp},e_{pq}\mid p,q\in Q_0\}$. 

The isotypical components of $V^m$ are $\Hom(E_{t(a)},E_{h(a)})$, 
$a\in Q_1$. Note that $m>m(Q)$ implies $m>\dim\Hom(E_p,E_q)$, and 
the semi-stability of the tautological bundles $\cal E_q$ follows 
from corollary \ref{cor:stab-bdl}.
\end{proof}

We wish to point out the following shortcoming: in this construction 
the length of the exceptional sequence equals the number of vertices 
of $Q$, which is independent of the multiplicity $m$. 
It follows that for large $m$ this sequence is certainly {\em not} 
a tilting object for $Y_{m}$ (compare this with remark \ref{rmk:length}).

%%===================================================================

\end{document}